\documentclass[11pt]{article}

\usepackage{array}
\usepackage{amsmath}
\usepackage{amssymb}
\usepackage{pstricks}
\usepackage{algorithm}
\usepackage{hyperref}
\usepackage[noend]{algorithmic}

\newcommand{\R}{\mathbb{R}}
\newcommand{\N}{\mathbb{N}}

\newcommand{\C}{\mathbb{C}}

\newcommand{\ac}{\`{a} }

\newtheorem{proposition}{Proposition}
\newtheorem{theorem}{Theorem}

\newtheorem{lemma}{Lemma}
\newtheorem{corollary}{Corollary}

\def\fdimo{~\hfill$\Box$\medbreak}

\begin{document}
	
\title{Local limit laws for symbol statistics\\
	in bicomponent rational models
\footnote{This report includes and improves the results presented in \cite{glv18,glv19a,glv19b}.}}
\author{Massimiliano Goldwurm$^{(1)}$, Jianyi Lin$^{(2)}$, Marco Vignati$^{(1)}$}

\date{ }

\maketitle

\begin{center}
{\small (1) Dipartimento di Matematica, Universit\ac degli Studi di Milano\\
	via Saldini 50, 20133 Milano, Italy \\
	 \{massimiliano.goldwurm,marco.vignati\}@unimi.it \\
 (2) Dipartimento di Scienze Statistiche, Universit\ac Cattolica del Sacro Cuore\\
	largo Gemelli 1, 20123 Milano, Italy \\
	jianyi.lin@unicatt.it}
\end{center}

\begin{abstract}
We study the local limit distribution of the number of occurrences of a symbol
in words of length $n$ generated at random in a regular language according to a rational stochastic model.
We present an analysis of the main local limits when the finite state automaton 
defining the stochastic model consists of two primitive components.
The limit distributions depend on several parameters and conditions, such as  
the main constants of mean value and variance of our statistics associated with the two components,
and the existence of communications from the first to the second component.
The convergence rate of these results is always of order $O(n^{-1/2})$.
We also prove an analogous $O(n^{-1/2})$ convergence rate to a Gaussian density
of the same statistic whenever the stochastic models only consists of one (primitive) component.
\end{abstract}

\noindent
{\bf Keywords}: limit distributions, local limit laws, pattern statistics, regular languages.

\section{Introduction}
\label{sec:introd}
In this work we present some local limit laws concerning the number of symbol occurrences in words
of given length chosen at random in a regular language under suitable probabilistic hypotheses.
Here the random words are assumed to be generated according to a \emph{rational stochastic model} as introduced in \cite{bcgl03},
which is defined by a (nondeterministic) finite state automaton with real positive weights on transitions.
In this setting the probability of generating a word $w$ is proportional to
the total weight of the transitions labelled by $w$ that are accepted by the automaton;
thus, the language recognized by the automaton is the family of all words having non-null probability to be generated.
This model is quite general, it includes as special cases
the traditional Bernoullian and Markovian sources \cite{nsf02,rs98}. 
It also includes the generation of random words of length $n$ in any regular language 
under uniform distribution:
this occurs when the automaton is unambiguous and all transitions have the same weight.

In order to fix ideas, consider a weighted finite state automaton $\cal A$
over an input alphabet including a special symbol $a$ and, for every $n\in\N$, let $Y_n$ be
the symbol statistics representing the number of occurrences of $a$ in a word of length $n$
generated at random according to the rational model defined by $\cal A$.
We are interested in the asymptotic properties of the sequence of random variables
$\{Y_n\}$.
These properties turn out to be of interest for the analysis of regular 
patterns occurring in words generated by Markovian models \cite{bcgl03,nsf02,rs98} and
for the asymptotic estimate of the coefficients of rational series in commutative variables \cite{bcgl03,bcgl06}.
They are also related to the descriptional complexity of languages and computational models \cite{bmr14},
as well as to the analysis of additive functions defined on regular languages \cite{gr07a}.
Clearly, the behaviour of $\{Y_n\}$ depends on automaton $\cal A$.
It is known that if $\cal A$ has a primitive transition matrix then $Y_n$ has a Gaussian limit distribution \cite{bcgl03,nsf02} and, under a suitable aperiodicity condition, it also satisfies a local limit theorem \cite{bcgl03}.
The limit distribution of $Y_n$ in the global sense is known also when the transition matrix of $\cal A$ consists of two primitive components \cite{dgl04}.

Here we improve the results of \cite{bcgl03,dgl04} presenting an analysis 
of the local limits of $\{Y_n\}$ when the transition matrix of $\cal A$ consists of two primitive components.
At the cost of adding suitable aperiodicity conditions, we prove that the main
convergence properties in distribution obtained in \cite{dgl04} also hold true in the local sense
with a convergence rate of the order $O(n^{-1/2})$.
We also obtain an analogous $O(n^{-1/2})$ rate of convergence for the Gaussian
local limit law of $\{Y_n\}$ when the rational model consists of only one primitive component,
so refining the result obtained in \cite{bcgl03}.
We recall that finding a tight convergence rate in central limit theorems is a natural goal 
of research \cite{hw98a}, which measures the approximation speed of the probability
values to the prescribed expression.

Our proofs are based on the analysis of the characteristic function of $Y_n$.
In particular we apply Laplace's method to evaluate the classical integral expression of this function that yields the probability values of interest.
This method is a general strategy to evaluate asymptotic expressions of
integrals depending on a growing parameter (see for instance \cite[Sec. B.6]{fs09}).
It is also used in the literature for different purposes, for instance to
prove local versions of the Central Limit Theorem \cite[Sec. 42]{gn97} or
for applications of the Saddle Point Method to combinatorial problems \cite[Ch. VIII]{fs09}.
The main difference with respect to these classical approaches is that
in our work Laplace's method is often applied to non-Gaussian integrals,
yielding (local) limit distributions that are not normal.

The material we present is organized as follows.
In the next section we define the problem and fix our notation.
In Section \ref{sec:primitive} we consider the rational stochastic models with primitive transition matrix and show how in this case Laplace's method leads to a local limit law of Gaussian type
for $\{Y_n\}$ with a convergence rate of the order $O(n^{-1/2})$.
In the same section we restate precisely an aperiodicity condition for the primitive stochastic model that is also necessary for the subsequent analysis.
In Section \ref{sec:bicomponente} we study the behaviour of $Y_n$ when the
rational stochastic model consists of two primitive components and has one or more transitions from the first to the second component.
In this case the following points summarize our results:
\begin{description}
\item[1.] If the two components have different main eigenvalues (and hence there is a dominant component), a Gaussian local limit law holds true.
Here the aperiodicity condition is assumed only for the dominant component.
\item[2.] If the two components have equal main eigenvalue (equipotent model) the 
limit distribution depends on the values of four constants: $\beta_1$,  $\gamma_1$ and $\beta_2$, $\gamma_2$, representing the leading terms of mean value and variance of our statistics associated with the first and the second component, respectively. In this case we assume the aperiodicity condition
for both components and the results are as follows:
\begin{description}
\item[2a.] If $\beta_1\neq \beta_2$ we get a local limit law for $Y_n/n$ towards a uniform density.
\item[2b.] If $\beta_1 = \beta_2$ but $\gamma_1 \neq \gamma_2$ we get a local limit law for 
$(Y_n - \beta n)/\sqrt{n}$
towards a suitable mixture of Gaussian densities.
\item[2c.] If $\beta_1 = \beta_2$ and $\gamma_1 = \gamma_2$ we obtain again
a Gaussian local limit.
\end{description}
\end{description}

In Section \ref{sec:sum} we study the behaviour of $Y_n$ in the
bicomponent models that have no communication between the components.
Also in this case, if there is a dominant component we get a Gaussian local limit.
On the contrary, in the equipotent case we get the following properties:
\begin{description}
	\item[3a.] If $\beta_1\neq \beta_2$ and/or $\gamma_1 \neq \gamma_2$ 
	we get a local limit law for $Y_n$ towards a convex linear combination of two Gaussian densities;
	\item[3b.] If $\beta_1 = \beta_2$ and $\gamma_1 = \gamma_2$ we obtain again
	a Gaussian local limit for $Y_n$.
\end{description}

All the local limit laws obtained above hold with a convergence rate $O(n^{-1/2})$.

Finally, in the last section we summarize and compare our results in a suitable table, also discussing possible goals for future investigations.

\section{Problem setting}
\label{sec:problemsetting}
As usual we denote by $\{a,b\}^*$ the set of all words over the binary alphabet $\{a,b\}$, 
including the empty word $\varepsilon$.
Together with the operation of concatenation between words, $\{a,b\}^*$ forms a monoid,
called \emph{free monoid} over $\{a,b\}$.
For every word $w\in \{a,b\}^*$ we denote by $|w|$ the length of $w$ 
and by $|w|_a$ the number of occurrences of $a$ in $w$.
For each $n\in\N$, we also represent by $\{a,b\}^n$ the set
$\{w\in \{a,b\}^*:|w|=n\}$. 
A \emph{formal series} in the non-commutative variables $a,b$ 
is a function $r:\{a,b\}^* \rightarrow \R_+$, where $\R_+ =\{x\in\R\mid x\geq 0\}$, 
and for every $w\in \{a,b\}^*$ we denote by $(r,w)$ the value of $r$ at $w$.
Such a series $r$ is called \emph{rational}
if for some integer $m>0$ there is a 
monoid morphism $\mu : \{a,b\}^* \rightarrow \R_+^{m\times m}$ and two (column) arrays
$\xi, \eta \in \R_+^m$, 
such that $(r,w) = \xi'\mu(w) \eta$, for every $w\in \{a,b\}^*$ \cite{br88,ss78}. 
In this case, as the morphism $\mu$ is generated by matrices
$A=\mu(a)$ and $B=\mu(b)$, we say that the 4-tuple
$(\xi,A,B,\eta)$ is a \emph{linear representation} of $r$ of size $m$.
Clearly, such a 4-tuple can be considered as a finite state automaton over the alphabet $\{a,b\}$,
with transitions (as well as initial and final states) weighted by positive real values.
Throughout this work we assume that the set $\{w\in\{a,b\}^n : (r,w)>0\}$ is not empty
for every $n\in\N_+$ (so that $\xi\neq 0\neq \eta$), and that
$A$ and $B$ are not null matrices, i.e. $A\neq [0] \neq B$.
Then we can consider the probability measure $\mbox{Pr}$ over the set $\{a,b\}^n$
given by
$$
\mbox{Pr}(w) = \frac{(r,w)}{\sum_{x\in\{a,b\}^n} (r,x)} =
\frac{\xi'\mu(w)\eta}{\xi'(A+B)^n\eta} \qquad 
\ \forall \ w\in \{a,b\}^n
$$
Note that, if $r$ is the characteristic series of a language $L\subseteq \{a,b\}^*$
then $\mbox{Pr}$ is the uniform probability function over the set
$L\cap \{a,b\}^n$.
Thus we can define the random variable (r.v.) $Y_n=|w|_a$, where $w$ is chosen at random in $\{a,b\}^n$
with probability $\mbox{Pr}(w)$.
As $A\neq [0] \neq B$, $Y_n$ is not a degenerate random variable.
It is clear that, for every $k\in\{0,1,\ldots,n\}$,
$$
p_n(k) := \mbox{Pr}(Y_n =k) = \frac{\sum_{|w|=n,|w|_a=k} (r,w)}{\sum_{w\in \{a,b\}^n} (r,w)}
$$
Since $r$ is rational also the previous probability can be expressed by using
its linear representation.
It turns out that
\begin{equation}
p_n(k) = \frac{[x^k]\xi' (Ax+B)^n \eta}{\xi' (A+B)^n \eta}
\qquad \forall \ k\in\{0,1,\ldots,n\}
\end{equation}
where, as usual, for any function $G$ analytic in a neighbourhood of $0$ with
series expansion $G(x) = \sum_{n=0}^{+\infty} g_nx^n$,
we denote by $[x^k] G(x)$ the coefficient $g_k$, for every $k\in \N$.

For sake of brevity we say that $Y_n$ is \emph{defined} by the linear representation
$(\xi,A,B,\eta)$.
To deal with the characteristic function $\Psi_n(t)$ of $Y_n$ we introduce the map $h_n(z)$ given by
\begin{equation}
\label{accaenne}
h_n(z) = \xi' (Ae^z + B)^n\eta \qquad \forall \ z \in\C 
\end{equation}
and hence we have
\begin{equation}
\label{fcaratteristica}
\Psi_n(t) = \sum_{k=0}^n p_n(k) e^{itk} = \frac{\xi' (Ae^{it} + B)^n\eta}{\xi' (A+B)^n\eta} = \frac{h_n(it)}{h_n(0)} \qquad \forall \ t \in\R
\end{equation}
As a consequence, mean value and variance of $Y_n$ can be evaluated by
\begin{equation}
\label{mediavar}
\mbox{E}(Y_n) = \frac{h_n'(0)}{h_n(0)} , \ \ 
\mbox{Var}(Y_n) = \frac{h_n''(0)}{h_n(0)} - \left( \frac{h_n'(0)}{h_n(0)} \right)^2
\end{equation}

Here we are mainly interested in the local limit properties of $\{Y_n\}$.
To compare the notion of local convergence with the traditional one,
we recall that a sequence of r.v.'s $\{X_n\}$ {\it converges in distribution} (or in law) 
to a random variable $X$ of distribution function $F$ if
$\ \lim_{n\rightarrow +\infty} \mbox{Pr}(X_n \leq x) = F(x)\ $,
for every $x\in \R$ of continuity for $F$ \cite{gn97}.
The central limit theorems yield classical examples of convergence in distribution 
to a Gaussian random variable. 

On the other hand, a local limit law for a sequence of discrete r.v.'s $\{X_n\}$ establishes,
as $n$ grows to infinity, an asymptotic expression for the probability values of $X_n$ 
depending an a given density function 
(see for instance \cite{be73,fs09,gn97}).
More precisely, assume that each $X_n$ takes value in $\{0,1,\ldots,n\}$.
We say that $\{X_n\}$ {\it satisfies a local limit law} of Gaussian type
if there are two real sequences $\{a_n\}$, $\{s_n\}$, 
where $\mbox{E}(X_n) \sim a_n$, $\mbox{Var}(X_n) \sim s_n^2$
and $s_n> 0$ for all $n$, 
such that for some real $\epsilon_n \rightarrow 0$, the relation
\begin{equation}
\label{gaussianlllaw}
\left| s_n \mbox{Pr}\left(X_n = k \right) \: - \: 
\frac{e^{-\left(\frac{k-a_n}{s_n}\right)^2/2}}{\sqrt{2\pi}} \right|
\: \leq \: \epsilon_n
\end{equation}
holds uniformly for every $k\in\{0,1,\ldots,n\}$ 
and every $n\in \N$ large enough.
Here, $\epsilon_n$ yields the {\it convergence rate} (or the speed) of the law.
A well-known example of such a property is given by the
de Moivre-Laplace local limit theorem \cite{gn97}.

Similar definitions can be given for other (non-Gaussian) 
types of local limit laws.
In this case the Gaussian density 
$e^{-x^2/2} / \sqrt{2\pi}$ appearing in (\ref{gaussianlllaw}) is replaced
by some density function $f(x)$; clearly, 
if $f(x)$ is not continuous at some points, the uniformity
with respect to $k$ must be adapted to the specific case.

We recall that in general convergence in distribution does not imply a local limit law;
usually, some further regularity condition is necessary to guarantee a local limit behaviour
\cite{fs09,gn97}.

\section{Primitive models}
\label{sec:primitive}

A relevant case occurs when $M=A+B$ is primitive, i.e. $M^k>0$ for some $k\in \N$ \cite{se81}.
In this case it is known that $Y_n$ has a Gaussian limit distribution 
and satisfies a local limit property \cite{bcgl03,nsf02}.
Here we improve this result, showing a rate of convergence $O(n^{-1/2})$;
we also recall some properties proved in \cite{bcgl03,bcgl06} 
that are useful in the following sections.

Since $M$ is primitive,
by Perron-Frobenius Theorem, it admits a real eigenvalue $\lambda > 0$ greater than the modulus of any other eigenvalue.
Moreover, strictly positive left and right eigenvectors $\zeta$, $\nu$ of $M$ 
with respect to $\lambda$
can be defined so that $\zeta'\nu=1$ \cite{se81}.
Thus, we can consider the function $u=u(z)$ implicitly defined by equation
$$ \mbox{Det}(Iu - Ae^z - B) = 0$$
subject to condition $u(0)= \lambda$.
It turns out that, in a neighbourhood of $z=0$, $u(z)$ is analytic,
it is a simple root of the characteristic polynomial of $Ae^z + B$, and
$|u(z)|$ is strictly greater than the modulus of all other eigenvalues of $Ae^z + B$.
Moreover, a precise relationship between $u(z)$ and function $h(z)$, defined in
(\ref{accaenne}), is proved in \cite{bcgl03} stating that
there are two positive constants $c$, $\rho$ and a function $r(z)$ analytic 
and non-null at $z=0$, such that
\begin{equation}
\label{sviluppo-h-eq}
h_n(z) = r(z)\; u(z)^n + O(\rho^n) \qquad \forall z \in \C : |z| \leq c
\end{equation}
where $\rho < |u(z)|$ and in particular $\rho < \lambda$.

Mean value and variance of $Y_n$ can be estimated from relations (\ref{sviluppo-h-eq})
and (\ref{mediavar}).
It turns out \cite{bcgl03} that the constants 
\begin{equation}
\label{betaegamma}
\alpha = \xi'\nu \zeta' \eta \ , \quad
\beta = \frac{u'(0)}{\lambda}  \quad \mbox{ and } \quad
\gamma = \frac{u''(0)}{\lambda} - \left(\frac{u'(0)}{\lambda} \right)^2
\end{equation}
are strictly positive and satisfy the equalities
$$
E(Y_n) = \beta n + O(1), \quad \beta = \frac{\nu' A \zeta}{\lambda}, \quad \mbox{ and } \quad
var(Y_n) = \gamma n + O(1)
$$

Other properties concern function $y(t) = u(it)/\lambda$,
defined for real $t$ in a neighbourhood of $0$.
In particular, there exists a constant $c>0$
for which identity (\ref{sviluppo-h-eq}) holds true,
satisfying the following relations \cite{bcgl03}:
\begin{equation}
\label{yt}
|y(t)| = 1 - \frac{\gamma}{2} t^2 + O(t^4), \  
\arg{y(t)} = \beta t + O(t^3), \   
|y(t)| \leq e^{-\frac{\gamma}{4}t^2} \qquad \forall\;|t| \leq c
\end{equation}
The behaviour of $y(t)$ can be estimated precisely when $t$ tends to $0$.
For any $q$ such that $1/3 < q < 1/2$  it can be proved \cite{bcgl03} that
\begin{equation}
\label{potenzayt}
y(t)^n = e^{-\frac{\gamma}{2}t^2n + i\beta t n} (1 + O(t^3)n)
\qquad \mbox{ for } |t| \leq n^{-q}
\end{equation}

Now, in order to prove a local limit property for $\{Y_n\}$ it is necessary to introduce an
aperiodicity assumption for the stochastic model, studied in more detail in \cite{bcgl06}.
To state this condition properly, 
consider the transition graph of the finite state automaton defined by matrices $A$ and $B$, i.e. 
the directed graph $G$ with vertex set $\{1,2,\ldots,m\}$ such that, for every
$i,j \in \{1,2,\ldots,m\}$, $G$ has an edge from $i$ to $j$ labelled by a letter
$a$ ($b$, respectively) whenever $A_{ij} > 0$ ($B_{ij} > 0$, resp.).
We can denote by $d$ the GCD of all differences in the number of occurrences of
$a$ in words representing (labels of) cycles of $G$ having equal length.
More formally, for every cycle $\cal C$ in $G$ let $\ell({\cal C})\in \{a,b\}^*$ be the word
obtained by concatenating the labels of all transitions in $\cal C$ in their order;
we define
$$ d = \mbox{GCD}
\{|\ell({\cal C}_1)|_a-|\ell({\cal C}_2)|_a : \  {\cal C}_1,{\cal C}_2 \mbox{ cycles in $G$ and } 
|{\cal C}_1| = |{\cal C}_2|\}
$$
Then, we say that the pair $(A,B)$ is \emph{aperiodic} if $d=1$.
Note that such a condition is often verified, for instance $d=1$ whenever $A_{ij} >0$ and $B_{ij} >0$
for two (possibly equal) indices $i,j$.
Moreover, it can be proved \cite{bcgl06} that $(A,B)$ is aperiodic if and only if,
for every real $t$ such that $0 < t < 2\pi$, we have
\begin{equation}
\label{condautovalori}
|\mu| < \lambda \ \ \mbox{ for every eigenvalue $\mu$ of } Ae^{it} + B
\end{equation}

\begin{theorem}
	\label{teo:locale_primitivo}
	Let $\{Y_n\}$ be defined by a linear representation $(\xi,A,B,\eta)$ such that
	matrix $M=A+B$ is primitive, $A \neq [0] \neq B$ and the pair $(A,B)$ is aperiodic.
	Moreover, let $\beta$ and $\gamma$ be defined by equalities (\ref{betaegamma}).
	Then, as $n$ tends to $+\infty$, the relation
	\begin{equation}
	\label{fine}
	\left| \sqrt{n} \mbox{Pr}\left(Y_n = k \right) \: - \: 
	\frac{e^{-\frac{(k-\beta n)^2}{2\gamma n}}}{\sqrt{2\pi\gamma}} \right|
	\: = \: \mbox{O}\left(n^{-1/2}\right)
	\end{equation}
	holds true uniformly for every $k \in \{0,1,\ldots,n\}$.
\end{theorem}
The statement is clearly meaningful when $k$, depending on $n$, varies so that
$x = \frac{k-\beta n}{\sqrt{2\gamma n}}$ lies in a finite interval.
In this case, we have
$\mbox{Pr}(Y_n = k) = \frac{e^{-x^2}}{\sqrt{2\pi\gamma n}} + O(n^{-1})$.

To prove the theorem, we study the characteristic function $\Psi_n(t)$ for $t \in [-\pi, \pi]$ 
by splitting this interval into three sets:
\begin{equation}
\label{intervalpr}
[-n^{-q},n^{-q}]\ , \qquad \{t\in\R : n^{-q} < |t| \leq c \}
\ , \qquad \{t\in\R : c < |t| \leq \pi\}\ ,
\end{equation}
where $c\in (0,\pi)$ is a constant satisfying relations (\ref{yt})
and  $q$ is an arbitrary value such that $\frac{1}{3} < q < \frac{1}{2}$.
We get the following three propositions where we always assume the hypotheses of
Theorem \ref{teo:locale_primitivo}.
\begin{proposition}
	\label{terzoint}
	For every $c \in(0, \pi)$ there exists $\varepsilon \in(0, 1)$ such that
	$$
	|\Psi_n(t)| = O(\varepsilon^n)\  \qquad  \forall\ t\in \R \ : \  c \leq |t| \leq \pi \ .
	$$
\end{proposition}
{\it Proof.}
First note that by property (\ref{condautovalori}), the aperiodicity of $(A,B)$ implies that, 
for every $c\in (0,\pi)$ there exists $\tau \in(0, \lambda)$
such that $|\mu| < \tau$ for every eigenvalue $\mu$ of $Ae^{it} + B$ 
and every $t\in\R$ satisfying $c \leq |t| \leq \pi$.
Also, by equality (\ref{accaenne}), the generating function of $\{h_n(it)\}_n$ is given by
$$
\sum_{n=0}^{+\infty} h_n(it) y^n \; = \; \xi' \left( I - (Ae^{it}+B) y \right)^{-1} \eta \; = \;
\frac{\xi' \mbox{Adj}\left( I - (Ae^{it}+B) y \right) \eta}{\mbox{Det}\left(I-(Ae^{it}+B) y\right)}
$$
and hence its singularities are the inverses of the eigenvalues of $Ae^{it} + B$.
As a consequence, $|h_n(it)| = O(\tau^n)$ whenever  $c \leq |t| \leq \pi$.
Moreover, from (\ref{sviluppo-h-eq}) we know that $h_n(0) = \Theta (\lambda^n)$
(\footnote{For any pair of sequences $\{f_n\}\subset \C$ and $\{g_n\}\subset \R_+$, we write 
	$f_n=\Theta(g_n)$ whenever there are two positive values $a$, $b$ such that 
	$a g_n \leq |f_n| \leq b g_n$ for every $n$ large enough.})
and hence, for some $\varepsilon \in (0,1)$, we have
$$
|\Psi_n(t)| = \left|\frac{h_n(it)}{h_n(0)} \right| = \frac{O(\tau^n)}{\Theta(\lambda^n)} = O(\varepsilon^n)
\qquad  \forall\ t\in \R \ : \  c \leq |t| \leq \pi \ 
$$
\fdimo
\begin{proposition}
	\label{secondoint}
	Let $c\in (0,\pi)$ satisfy relation (\ref{yt}).
	Then, for every $t\in \R$ such that $n^{-q} \leq |t| \leq c$ we have
	$$
	|\Psi_n(t)| = O\left( e^{-\frac{\gamma}{4}n^{1-2q}} \right)
	$$
\end{proposition}
{\it Proof.}
By relation (\ref{sviluppo-h-eq}), since  $y(0) = 1$, 
there exists $\rho \in (0, \lambda) $ such that
\begin{equation}
\label{psiditinzero}
\Psi_n(t) = \frac{h_n(it)}{h_n(0)} = 
\frac{r(it) \lambda^n y(t)^n + O(\rho^n)}{r(0) \lambda^n + O(\rho^n)}
\ \qquad \forall\ t\in \R \ : \ |t| \leq c 
\end{equation}
Since $r(z)$ is analytic in a neighbourhood of $0$, we have
$$
|\Psi_n(t)| = (1 + O(t)) |y(t)|^n + O(\varepsilon^n) \ ,\qquad 
\mbox{ for }\ 0<\varepsilon< 1 \ .
$$
Also, by inequality (\ref{yt}), we know that 
$|y(t)|^n \leq e^{-\frac{\gamma}{4}t^2n}$ whenever $|t| \leq c$.
Thus, the result follows by replacing this bound in the previous equation and recalling that 
$n^{-q} \leq |t| \leq c$.

\fdimo

\begin{proposition}
	\label{intpsint}
	For any $q$ such that $1/3 < q < 1/2$, we have
	$$
	\int_{|t|\leq n^{-q}}
	\left| \Psi_n(t) - e^{-\frac{\gamma}{2}t^2n +i\beta tn}  \right| dt \: = \: O(n^{-1})
	$$
\end{proposition}
{\it Proof.}
Reasoning as in Proposition \ref{secondoint},
from relation (\ref{psiditinzero}) we know that
$\Psi_n(t) = (1 + O(t)) y(t)^n + O(\varepsilon^n)$
for some $\varepsilon \in (0,1)$ and every $t\in \R$ such that $|t| \leq c$.
Thus, applying relation (\ref{potenzayt}) and recalling that $nO(t^3) = o(1)$ for $|t| \leq n^{-q}$, we get
$$
\Psi_n(t) = (1 + O(t) + nO(t^3)) e^{-\frac{\gamma}{2}t^2n +i\beta tn} +
O(\varepsilon^n)
\  \qquad \forall \ t\in \R \ : \ |t| \leq n^{-q}
$$
Thus, computing directly the primitives of simple functions, we obtain
\begin{eqnarray*}
\nonumber
\lefteqn{\int_{|t|\leq n^{-q}}
	\left| \Psi_n(t) -  e^{-\frac{\gamma}{2}t^2n +i\beta tn}  \right| dt = } \\
\nonumber & = &
\int_{|t|\leq n^{-q}} |O(t) + nO(t^3)|\; e^{-\frac{\gamma}{2}t^2n} dt +
O(\varepsilon^n)  \ = \ O\left( n^{-1} \right) \qquad \qquad \qquad
~\hfill \Box
\end{eqnarray*}

Now, we are able to prove the result of this section.

\noindent
{\bf Proof of Theorem \ref{teo:locale_primitivo}.}
It is well-known \cite{gn97} that $p_n(k)= \mbox{Pr}\left\{Y_n = k \right\}$ 
can be computed from the inversion formula
\begin{equation}
\label{inversionform}
p_n(k) \: = \: \frac{1}{2\pi} \int_{-\pi}^{\pi} \Psi_n(t) e^{-itk} dt
\qquad \forall k\in\{0,1,\ldots,n\}
\end{equation}
To evaluate that integral, let us split $[-\pi,\pi]$ into the three sets 
defined in (\ref{intervalpr}).
Then, by Propositions \ref{terzoint} and \ref{secondoint}, 
for some $\varepsilon \in (0,1)$ we obtain
\begin{equation}
\label{pnkint}
p_n(k) \: = \: \frac{1}{2\pi} \int_{|t|\leq n^{-q}} \Psi_n(t) e^{-itk} dt +
O\left( e^{-\frac{\gamma}{4}n^{1-2q}} \right) + 
O(\varepsilon^n) 
\end{equation}
Moreover, setting $v (=v_{k,n}) = \frac{k-\beta n}{\sqrt{\gamma n}}$,
by Proposition \ref{intpsint} we have
\begin{equation}
\label{due}
\int_{|t|\leq n^{-q}} \Psi_n(t) e^{-itk} dt \ = \ 
\int_{|t|\leq n^{-q}}  e^{-\frac{\gamma}{2}t^2n - i t v \sqrt{\gamma n}} dt + O(n^{-1})
\end{equation}
By a standard computation the first term on the right-hand side becomes 
\begin{eqnarray}
\nonumber
\int_{|t|\leq n^{-q}}  e^{-\frac{\gamma}{2}t^2n - i t v \sqrt{\gamma n}} dt
& = & \frac{1}{\sqrt{\gamma n}} \left( \int_{-\infty}^{+\infty} e^{-\frac{x^2}{2}-ivx} dx -
\int_{|x| > n^{\frac{1}{2}-q} \sqrt{\gamma}} e^{-\frac{x^2}{2}-ivx} dx  \right) \\
\label{uno}
& = &  \frac{1}{\sqrt{\gamma n}} \left( \sqrt{2\pi}\  e^{-\frac{v^2}{2}}
+ O(e^{-\frac{\gamma}{2}n^{1-2q} } ) \right)
\end{eqnarray}
where one recognizes in the second integral the characteristic function of a 
Gaussian random variable.
Thus, the result follows by replacing (\ref{uno}) in (\ref{due}) and
(\ref{due}) in (\ref{pnkint}).
\fdimo

\section{Bicomponent models}
\label{sec:bicomponente}

In this section we study the behaviour of $\{Y_n\}_{n\in\N}$ defined by a linear representation
$(\xi,A,B,\eta)$ of size $m$ such that the matrix $M=A+B$ consists of two irreducible components.
Formally, there are two linear representations, $(\xi_1,A_1,B_1,\eta_1)$ and $(\xi_2,A_2,B_2,\eta_2)$, of size 
$m_1$ and $m_2$ respectively, where $m=m_1+m_2$, such that:

\begin{enumerate}
\label{ipotesibicomp}
\item For some $A_0, B_0 \in \R_+^{m_1\times m_2}$ we have
	$$
	\label{2rappr} \xi' = (\xi'_1, \xi'_2) , \quad 
	A = \left( \begin{array}{cc} A_1 & A_0 \\ 0 & A_2 \end{array} \right), \quad 
	B = \left( \begin{array}{cc} B_1 & B_0 \\ 0 & B_2 \end{array} \right), \quad
	\eta=\left( \begin{array}{c} \eta_1 \\ \eta_2 \end{array} \right)
	$$ 
\item $M_1 =A_1 + B_1$ and $M_2 =A_2 + B_2$ are irreducible matrices and we denote by
$\lambda_1$ and $\lambda_2$ the corresponding Perron-Frobenius eigenvalues;
\item $\xi_1 \neq 0 \neq \eta_2$ and
matrix $M_0 = A_0 + B_0$ is different from $[0]$.
\end{enumerate}

Note that condition 2 is weaker than a primitivity assumption for $M_1$ and $M_2$.
Moreover, condition 3 avoids trivial situations and guarantees a (non-vanishing) communication from
the first to the second component.

A typical example of a formal series $r$ with a linear representation of this kind is given by the product
of two rational formal series $r_1$, $r_2$, both having an irreducible linear representation,
i.e. $r = r_1 \cdot r_2$, meaning that $(r,w) = \sum_{w=xy} (r_1,x)(r_2,y)$ for every $w\in\{a,b\}^*$.

Under these hypotheses the limit distribution of $\{Y_n\}$ first depends on whether
$\lambda_1 \neq \lambda_2$ or $\lambda_1=\lambda_2$.
If $\lambda_1 \neq \lambda_2$ there is a dominant component, 
corresponding to the maximum between $\lambda_1$ and $\lambda_2$, which determines the asymptotic behaviour of $\{Y_n\}$.
If $\lambda_1=\lambda_2$ the two components are equipotent and they both contribute to
the limit behaviour of $\{Y_n\}$.
In both cases the corresponding characteristic function has some common properties
we now recall briefly.

For $j=1,2$, let us define
$h_n^{(j)}(z)$, $u_j(z)$, $y_j(t)$, $\beta_j$, and $\gamma_j$, respectively, as the values 
$h_n(z)$, $u(z)$, $y(t)$, $\beta$, $\gamma$ referred to component $j$.
We also define $H(x,y)$ as the matrix-valued function given by
\begin{eqnarray}
\nonumber
H(x,y) & = & \sum_{n=0}^{+\infty} (Ax+B)^n y^n= 
\left[ \begin{array}{cc} H^{(1)}(x,y) & G(x,y) \\ 
0 & H^{(2)}(x,y) \end{array} \right]\ , \quad \mbox{ where} \\
\label{H1eH2}
H^{(1)}(x,y) & = & 
\frac{\mbox{Adj}\left(I-(A_1x+B_1)y\right)}{\mbox{Det}\left(I-(A_1x+B_1) y\right)}, \ 
H^{(2)}(x,y) \ = \ 
\frac{\mbox{Adj}\left(I-(A_2x+B_2)y\right)}{\mbox{Det}\left(I-(A_2x+B_2) y\right)}, \\ 
\nonumber 
\mbox{and} \ & & 
G(x,y) \ = \ H^{(1)}(x,y) \; (A_0 x + B_0)y \; H^{(2)}(x,y) \ .
\end{eqnarray}
Thus, the generating function of $\{h_n(z)\}_n$ satisfies the following identities
\begin{equation}
\label{fungenhnbicomp}
\sum_{n=0}^\infty h_n(z) y^n = \xi' H(e^z,y) \eta = \xi_1' H^{(1)}(e^z,y)  \eta_1  +
\xi_1' G (e^z,y)  \eta_2  +  \xi_2' H^{(2)} (e^z,y) \eta_2 
\end{equation}
Hence, setting $g_n(z)  = [y^n]\xi_1' G(e^z,y) \eta_2\ $,
we obtain
\begin{equation}
\label{acca_decomp} h_n(z) = h_n^{(1)}(z) + g_n(z) + h_n^{(2)}(z)
\end{equation}

\subsection{Dominant case}

Under the previous hypotheses let us further assume  
$\lambda_1 > \lambda_2$ and $M_1$ aperiodic (and hence primitive).
In this case we say that $\{Y_n\}$ is defined 
in a \emph{dominant communicating bicomponent model} with $\lambda_1 > \lambda_2$.
Under these assumptions it is known that if $A_1\neq 0 \neq B_1$ then $0<\beta_1 < 1$, $0<\gamma_1$ and 
$\frac{Y_n-\beta_1 n}{\sqrt{\gamma_1 n}}$ converges in distribution to a normal r.v.
of mean value $0$ and variance $1$ \cite{dgl04}.
Here we show a Gaussian local limit law for $\{Y_n\}$
with a convergence rate $O(n^{-1/2})$,
under the further hypothesis that $(A_1,B_1)$ is aperiodic.
Note that, by definition, the aperiodicity of $(A_1,B_1)$ implies
$A_1 \neq 0 \neq B_1$ (and hence $0<\beta_1 < 1$, $0<\gamma_1$).
The proof is similar to that one of Theorem \ref{teo:locale_primitivo}
and here we present a brief outline.

\begin{theorem}
	\label{teo:locale_bicompdom}
	Let $\{Y_n\}$ be defined in a dominant communicating bicomponent model 
	with $\lambda_1 > \lambda_2$ and assume $(A_1,B_1)$ aperiodic.
	Then, as $n$ tends to $+\infty$, the relation 
	$$
	\left| \sqrt{n} \mbox{Pr}\left(Y_n = k \right) \: - \: 
	\frac{e^{-\frac{(k-\beta_1 n)^2}{2\gamma_1 n}}}{\sqrt{2\pi\gamma_1}} \right|
	\: = \: \mbox{O}\left(n^{-1/2}\right)
	$$
	holds true uniformly for every $k\in \{0,1,\ldots,n\}$.
\end{theorem}
Also in this case the statement is significant when 
$ \frac{(k-\beta_1 n)^2}{2\gamma_1 n} $ remains in a finite interval
\footnote{Even if not explicitly mentioned, similar observation also holds for the other local limit theorems presented in this work, whenever the limit function is continuous.}.

As in the previous section, the proof is based on the analysis of the characteristic function $\Psi_n(t)$ for $t$ lying in the three sets
$\ |t| \leq n^{-q}\ $, $\ n^{-q} < |t| \leq c\ $ and $\ c < |t| \leq \pi\ $,
where $c\in (0,\pi)$ is a suitable constant and $q$ is an arbitrary real value such that $\frac{1}{3} < q < \frac{1}{2}$.
First, let us consider the third set.
\begin{proposition}
	\label{terzointbicomp} 
	Under the hypotheses of Theorem \ref{teo:locale_bicompdom}, 
	for every $c\in (0,\pi)$ there exists $\varepsilon \in (0,1)$ such that
	$$
	|\Psi_n(t)| = O(\varepsilon^n)\ , \qquad \forall\ t\in \R \ : \  c \leq |t| \leq \pi \ .
	$$
\end{proposition}
{\it Proof.}
Note that, by relations (\ref{H1eH2}) and (\ref{fungenhnbicomp}), for every $t \in \R$ the singularities of the generating function
$\xi' H(e^{it},y)\eta = \sum_{n=0}^\infty h_n(it) y^n$
are the inverses of the eigenvalues of $A_1e^{it} + B_1$ and  $A_2e^{it} + B_2$.
Assuming $c \leq |t| \leq \pi$, the first ones are in modulus smaller than $\lambda_1$ by condition (\ref{condautovalori}), 
while the second ones are in modulus smaller or equal to $\lambda_2$ as a consequence of
Perron-Frobenius Theorem for irreducible matrices \cite[Ex. 1.9]{se81}.
Thus, since $\lambda_1 > \lambda _2$, for some positive $\tau < \lambda_1$ we have
$|h_n(it)| = O(\tau^n)$ for all real $t$ such that $c \leq |t| \leq \pi$.
For the same argument it is clear that
$h_n(0) = \Theta(\lambda_1^n)$, and hence
the result follows by reasoning as in the proof of Proposition \ref{terzoint}.
\fdimo

Concerning the other two subsets, 
i.e. $|t| \leq n^{-q}$ and $n^{-q} < |t| \leq c$, 
one can study the behaviour of $h_n(z)$ in a neighbourhood of $z=0$.
Reasoning as in the previous section it is easy to show that there are
two positive constants $c$, $\rho$ such that
$$
h_n(z) = r(z)\; u_1(z)^n + O(\rho^n) \qquad \forall z \in \C : |z| \leq c
$$
where $\rho < |u_1(z)|$ and $r(z)$ is a function analytic and non-null for $z\leq c$.
In particular $\rho < \lambda_1$ and $h_n(0) = r(0)\lambda_1^n + O(\rho^n)$.
These properties allow us to argue as in Section \ref{sec:primitive}, 
replacing the values $\lambda$, $\beta$, $\gamma$, $y(t)$ respectively by
$\lambda_1$, $\beta_1$, $\gamma_1$ and $y_1(t)$, thus proving
two statements equivalent to Propositions \ref{secondoint} and \ref{intpsint}, respectively.
The proof of Theorem \ref{teo:locale_primitivo} can be modified in the same way
and this concludes the proof of Theorem \ref{teo:locale_bicompdom}.

\subsection{Equipotent case}
\label{subsec:equipotent}

Now let us consider the equipotent case.
Formally, let $\{Y_n\}$ be defined by a linear representation $(\xi,A,B,\eta)$ satisfying
conditions 1, 2, 3 above,
assume $\lambda_1=\lambda_2 = \lambda$ and let
both matrices $M_1$, $M_2$ be aperiodic (and hence primitive).
Under these hypotheses we say that $\{Y_n\}$ is defined in an \emph{equipotent communicating bicomponent model}.
The limit distribution of $\{Y_n\}$ in this case is studied in \cite{dgl04} and
depends on the parameters $\beta_1$, $\beta_2$, $\gamma_1$, $\gamma_2$.
Here we extend those results to local limit properties with a suitable convergence rate
under the further assumption that both pairs $(A_1,B_1)$ and $(A_2,B_2)$ are aperiodic
(again, such a hypothesis guarantees that $0 < \beta_j < 1$ and $0 < \gamma_j $, for both $j=1,2$).
To this end we first recall some useful properties of the characteristic function of $Y_n$
under these hypotheses \cite{dgl04}.

First, consider equality (\ref{acca_decomp}) and note that both $h_n^{(1)}(z)$ and 
$h_n^{(2)}(z)$ satisfy (an analogue of) relation (\ref{sviluppo-h-eq}).
Moreover, to evaluate $\{g_n(z)\}$ observe that its generating function
$\xi_1' G(e^z,y) \eta_2$, for every complex $z$ in a neighbourhood of $0$, 
has the singularities of smallest modulus at points
$y=u_1(z)^{-1}$ and $y=u_2(z)^{-1}$.
Thus, for a suitable $c>0$ we can write
$$
\xi_1' G(e^z,y) \eta_2 = \frac{s(z)y}{(1-u_1(z)y)(1-u_2(z)y)} + L(z,y)
\qquad \forall z \in \C : |z| \leq c
$$
where again $s(z)$ is a function analytic and non-null for $z\leq c$,
and $L(z,y)$ only admits singularities of modulus strictly greater than 
$|u_1(z)|^{-1}$ and $|u_2(z)|^{-1}$.
This implies
\begin{equation}
\label{gennedizeta}
g_n(z) = s(z) \sum_{j=0}^{n-1} u_1(z)^j u_2(z)^{n-1-j} + O(\rho^n)
\qquad \forall z \in \C : |z| \leq c
\end{equation}
where $\rho < \max\{|u_1(z)|, |u_2(z)|\}$.
Replacing this expression in (\ref{acca_decomp}), one gets
\begin{equation}
\label{svihnbicomp}
h_n(z) = s(z) \sum_{j=0}^{n-1} u_1(z)^j u_2(z)^{n-1-j} + 
O(u_1(z)^n) + O(u_2(z)^n)
\qquad \forall z \in \C : |z| \leq c
\end{equation}
This equality has two consequences. 
First, since $u_1(0) = \lambda = u_2(0)$, it implies
\begin{equation}
\label{accan0bicomp}
h_n(0) = s(0) n \lambda ^{n-1}(1 +  O(1/n)) \qquad \mbox{($s(0)\neq 0$)}
\end{equation}
Second, if $u_1(z) \neq u_2(z)$ for some $z\in\C$ satisfying $ 0<|z|\leq c$,
one gets
\begin{equation}
\label{accanzetabicomp}
h_n(z) = s(z)\frac{u_1(z)^n - u_2(z)^n}{u_1(z) - u_2(z)}
+ O(u_1(z)^n) + O(u_2(z)^n)
\end{equation}

Finally, assuming the above aperiodicity condition and 
reasoning as in Proposition \ref{terzoint}, 
the following property can be proved by using 
relations (\ref{H1eH2}) and (\ref{accan0bicomp}) .
\begin{proposition}
	\label{terzointbic}
	Let $\{Y_n\}$ be defined in an equipotent communicating bicomponent model
	and let both pairs $(A_1,B_1)$ and $(A_2,B_2)$ be aperiodic.
	Then, for every $c \in (0,\pi)$ there exists $\varepsilon \in (0,1)$ such that
	$\ |\Psi_n(t)| = O(\varepsilon^n)\ $ 
	for all $t\in\R$ satisfying $c \leq |t| \leq \pi$.
\end{proposition}

\subsubsection{Local limit with different $\beta$'s}

In this subsection we assume an equipotent communicating bicomponent model
with $\beta_1 \neq \beta_2$.
In this case it is known that $Y_n/n$ converges in distribution to a uniform r.v. over the interval of extremes $\beta_1$, $\beta_2$.
Here we prove a local limit theorem toward the corresponding density function
with a convergence rate of order $O(n^{-1/2})$.
To this end, in view of Proposition \ref{terzointbic}, we study the characteristic function 
$\Psi_n(t)$ for $|t| \leq c$,
where $c\in (0,\pi)$ is a constant for which identity (\ref{accanzetabicomp}) holds true.
Recall that in such a set both functions
$y_1(t)= u_1(it)/\lambda$ and $y_2(t)= u_2(it)/\lambda$ satisfy relations
(\ref{yt}),
and hence for every real $t$ such that $|t| \leq c$, we have
\begin{gather}
\label{ytbicomp}
y_j(t) = 1 + i\beta_jt + O(t^2) \ , \qquad j=1,2 \\
\label{ytbound}
|y_j(t)| \leq e^{-\frac{\gamma_j}{4}t^2} \ , \qquad j=1,2
\end{gather}
As a consequence, we may assume the following relation for every $t\in\R$ such that $0<|t| \leq c$:
\begin{equation}
\label{psinbicomp}
\Psi_n(t) = \frac{h_n(it)}{h_n(0)} =  
\frac{1+O(t)}{1+O(1/n)} 
\left( \frac{y_1(t)^n - y_2(t)^n}{i\;(\beta_1-\beta_2)\;t n} \right)
+ \sum_{j=1,2} O\left( \frac{y_j(t)^n}{n} \right) 
\end{equation}
Now, for such a constant $c$, let us split the interval $[-c,c]$ 
into sets $S_n$ and $V_n$ given by
{\small
\begin{equation} 
\label{intervalli}
S_n = \left\{t\in\R : |t| \leq \frac{\log n}{\sqrt{n}}\, \right\},\; 
V_n = \left\{ t\in\R : \frac{\log n}{\sqrt{n}} < |t| \leq c \right\} 
\end{equation} } 
The behaviour of $\Psi_n(t)$ in $V_n$ is given by the following proposition, where
we assume an equipotent communicating bicomponent model with $\beta_1 \neq \beta_2$.

\begin{proposition}
	\label{secondointbic}
	It turns out that
	$\ |\Psi_n(t)| = o\left( n^{-3/2} \right)\ $
	for all $t \in V_n$.
\end{proposition}
{\it Proof.}
From equation (\ref{psinbicomp}), for every $t\in V_n$, we obtain  
$$
|\Psi_n(t)| \leq \frac{|y_1(t)|^n+|y_2(t)|^n}{\Theta(t n)} + \sum_{j=1,2} 
O\left( \frac{|y_j(t)|^n}{n} \right)  =
\sum_{j=1,2} o\left( \frac{|y_j(t)|^n}{\sqrt{n}} \right)
$$
Taking $a= \frac{\min\{\gamma_1,\gamma_2\}}{4}$ by relations (\ref{ytbound})
we get 
$
|\Psi_n(t)| = 
o\left( n^{-\frac{1}{2}-a(\log n)^2} \right)
$,
which proves the result.
\fdimo

Now, let us evaluate $\Psi_n(t)$ for $t\in S_n$.
To this end we need the following
\begin{lemma}
	\label{segnialt}
	For $k,m \in \N$, $k < m$, let $g:[2k\pi,2m\pi]\rightarrow \R_+$ 
	be a monotone function, and let
	$I_{k,m} = \int_{2k\pi}^{2m\pi} g(x) \sin x\; dx $. 
	Then:
	
	a) if $g$ is non-increasing we have 
	$\ 	0 \leq I_{k,m} \leq 2[g(2k\pi)-g(2m\pi)]$;
	
	b) if $g$ is non-decreasing we have 
	$\ 2[g(2k\pi)-g(2m\pi)] \leq I_{k,m} \leq 0 $.
	
	In both cases $|I_{k,m}| \leq 2|g(2k\pi)-g(2m\pi)|$.
\end{lemma}
{\it Proof.}
If $g$ is non-increasing, for each integer $j\in [k,m)$ we have $0 \leq I_{j,j+1} $ and
$$
I_{j,j+1} = \int_{2j\pi}^{(2j+1)\pi} g(x) \sin x \; dx - 
\int_{(2j+1)\pi}^{2(j+1)\pi} g(x) |\sin x| \; dx \leq
2[g(2j\pi) - g(2(j+1)\pi)]
$$
Thus a) follows by summing the expressions above for $j=k,\ldots,m-1$.
Part  b) is proved by applying a) to function
$h(x) = g(2m\pi) - g(x)$.
\fdimo

Going back to the analysis of $\Psi_n(t)$ in $S_n$,
let us define 
\begin{equation}
\label{cappaennediti}
K_n(t) = \frac{e^{-\frac{\gamma_1}{2}t^2n + i\beta_1 t n} - 
	e^{-\frac{\gamma_2}{2}t^2n + i\beta_2 t n}}{i(\beta_1-\beta_2) t n}
\end{equation}
and consider relation (\ref{psinbicomp}).
Since for $t\in S_n$ one has $nO(t^3) =o(1)$, relation (\ref{potenzayt})
applies to both $y_1(t)$ and $y_2(t)$, yielding
$$
y_j(t)^n  =  e^{-\frac{\gamma_j}{2}t^2n + i\beta_j t n} (1 + nO(t^3))
\qquad \forall \ t \in S_n , \quad j=1,2
$$
Replacing these values in (\ref{psinbicomp}), for some $a > 0$ one gets
\begin{equation}
\label{intermedio}
\Psi_n(t) \ = \ 
\left[ 1  + O(t) + n O(t^3) + O(1/n)  \right] K_n(t) + O(n^{-1} e^{-at^2 n})
\qquad \forall\ t\in S_n
\end{equation}
Such an equality allows to determine the properties of $\Psi_n(t)$ in $S_n$.

\begin{proposition}
\label{primointbic}
Assume an equipotent communicating bicomponent model with $\beta_1 \neq \beta_2$ and let
$S_n$ and $K_n(t)$ be defined as in {\rm (\ref{intervalli})} and {\rm (\ref{cappaennediti})},
respectively.  Then, we have
$$
\left| \int_{S_n} \left( \Psi_n(t) - K_n(t) \right) dt \right|  \ = \ 
O\left(n^{-3/2}\right)
$$
\end{proposition}
{\it Proof.}
Integrating both sides of (\ref{intermedio}), we obtain
\begin{equation}
\label{intdiff}
\int_{S_n} \left[ \Psi_n(t) - K_n(t) \right] dt \ = \ 
\int_{S_n} \left\{ \left[ O(t) + nO(t^3) + O(n^{-1}) \right] K_n(t) + O(n^{-1} e^{-at^2 n})
\right\}  dt  
\end{equation}
In order to evaluate the integral in the right hand side
observe that, for any constant $a>0$ and $\tau_n = n^{-1/2}(\log n)$, we have
\begin{equation}
\label{baseuno}
\int_{S_n} e^{-at^2 n} dt \ \leq \ \frac{2}{\sqrt{n}} + 
2 \int_{n^{-1/2}}^{\tau_n} \sqrt{n} t\; e^{-at^2 n}\; dt \ = \ \Theta\left(n^{-1/2}\right)
\end{equation}
which implies (for a suitable $a>0$)
\begin{equation}
\label{basedue}
\left| \int_{S_n} t K_n(t) dt \right| \ \leq \ 
O \left( \int_{S_n} n^{-1} \; e^{-at^2 n} \; dt \right) = O\left( n^{-3/2} \right)
\end{equation}
Moreover, using similar bounds one gets 
\begin{eqnarray}
\nonumber
	\left| \int_{S_n} nt^3 K_n(t) dt \right| & = & 
	O\left( \int_{S_n} t^2 \; e^{-at^2 n} \; dt \right) = 
	O\left( n^{-3/2} + \int_{n^{-1/2}}^{\tau_n} t^2 \; e^{-at^2 n} \; dt \right) \\
	& = & \hspace{-0.2cm}
	O\left( \hspace{-0.1cm} n^{-3/2} + \frac{1}{n} \left\{ \left. -te^{-at^2 n} \right|_{n^{-1/2}}^{\tau_n} + \hspace{-0.1cm}
	 \int_{n^{-1/2}}^{\tau_n} e^{-at^2 n} \; dt \right\} \hspace{-0.1cm} \right) = 
	O\left( n^{-3/2} \right) 
	\label{basetre}
\end{eqnarray}
Using (\ref{baseuno}), (\ref{basedue}) and (\ref{basetre}) in (\ref{intdiff})
one easily see that the result is proved once we show
\begin{equation}
\label{intcappaennediti}
\int_{S_n} K_n(t) dt = O(1/n)
\end{equation}
To this end, define $\delta = \beta_1 - \beta_2$.
Since $\cos x$ and $\sin x$ are respectively even and odd function, we can write
\begin{equation*}
	\left| \int_{S_n} K_n(t) dt \right| \leq  
	\sum_{j=1,2} \left| \int_{S_n} \frac{e^{-\frac{\gamma_j}{2}t^2n + i\beta_j t n}-1}{i\delta tn} dt \right| = 
	\sum_{j=1,2} \left| \int_{S_n} \frac{e^{-\frac{\gamma_j}{2}t^2n}\sin(\beta_j tn)}{\delta tn} dt \right| 
\end{equation*}
Setting $u = \beta_j n t$, each integral in the last sum becomes (for some $a,b >0$)
\begin{equation}
\label{passocappa}
\frac{2}{|\delta|n} \int_0^{b n \tau_n} \frac{e^{-a\frac{u^2}{n}}\; \sin u}{u}\; du 
\ \leq \ 
\frac{2}{|\delta|n} \left\{ \int_0^{2\pi} \frac{sin u}{u}\; du + 
\int_{2\pi}^{b n \tau_n} \frac{e^{-a\frac{u^2}{n}}\; \sin u}{u}\; du \right\}
\end{equation}
Thus, we can apply Lemma \ref{segnialt} to $g(u) = u^{-1} e^{-a\frac{u^2}{n}}$
in the last expression, and get
\begin{equation*}
\int_{2\pi}^{bn \tau_n} u^{-1} e^{-a\frac{u^2}{n}}\; \sin u \; du \ = \ 
	2 \left( g(2\pi) - g(b\tau_n n) + o(1) \right) \ = \ \pi^{-1} + o(1)
\end{equation*}
Replacing this value in (\ref{passocappa}) we obtain equality (\ref{intcappaennediti})
and the proof is complete.
\fdimo

Now, we are able to prove the local limit in the present case.
Set $b_1 = \min\{\beta_1,\beta_2\}$, 
$b_2 = \max\{\beta_1,\beta_2\}$
and denote by $f_U(x)$ the density function of a uniform r.v. $U$ in the interval $[b_1,b_2]$,
that is
$$f_U(x) = \frac{1}{b_2-b_1} \chi_{[b_1,b_2]}(x)\ \qquad \forall x \in \R$$
where $\chi_I$ denotes the indicator function of interval $I\subset \R$.

\begin{theorem}
	\label{teo:bicompequipbetadiv}
	Let $\{Y_n\}_{n\in\N}$ be defined in an equipotent communicating bicomponent model
	with $\beta_1\neq \beta_2$ and assume aperiodic both pairs $(A_1,B_1)$ and $(A_2,B_2)$.
	Then, for $n$ tending to $+\infty$, $Y_n$ satisfies the relation
	\begin{equation}
		\label{converglocale}
		\left| n\; \mbox{Pr}(Y_n= k) \; - \; f_U(x) \right| \ = \ 
		O\left( n^{-1/2} \right)
	\end{equation}
	for every $k=k(n)\in\N$ such that $\lim_{n\rightarrow\infty}k/n = x$, 
	where $x$ is a constant different from $\beta_1$ and $\beta_2$.
\end{theorem}
{\it Proof.}
We start again from the inversion formula (\ref{inversionform}).
To evaluate the integral therein we split the interval $[-\pi,\pi]$ into the three sets 
$\{t\in\R : c < |t| \leq \pi\}$, $V_n$, and $S_n$, 
where $c>0$ is a constant for which relation (\ref{psinbicomp}) holds true, while
$S_n$ and $V_n$ are defined in equations (\ref{intervalli}).
Then, by Propositions \ref{terzointbic}, \ref{secondointbic}, \ref{primointbic}, we obtain
\begin{equation}
	\label{pnbic}
	p_n(k) = \frac{1}{2\pi} \int_{S_n}
	\left( \frac{e^{-\frac{\gamma_2}{2}t^2n + i\beta_2 t n} - 
		e^{-\frac{\gamma_1}{2}t^2n + i\beta_1 t n}}{it\;(\beta_2-\beta_1)\; n} \right) 
	e^{-ikt} dt + O\left( n^{-3/2}\right)
\end{equation}
Now, set $v=k/n$ and note that for $n\rightarrow +\infty$, $v$ converges to a constant $x$
different from $\beta_1$ and $\beta_2$.
Thus, defining
$$
\Delta_n(v) = \int_{S_n} 
\frac{e^{i(\beta_2 - v)tn - \frac{\gamma_2}{2} t^2n} - e^{i(\beta_1 - v)tn - 
		\frac{\gamma_1}{2} t^2n}}{i(\beta_2 - \beta_1) t} \; dt
$$
we are done once we prove that
\begin{equation}
	\label{finale}
	\Delta_n(v) = 2\pi f_U(x) + O(n^{-1/2})
\end{equation}
To this end, without loss of generality assume $\beta_1 < \beta_2 $
and set $\delta = \beta_2 - \beta_1$.
Then, $\Delta_n(v)$ is an integral of the difference between two functions of the form
$$
A_n(t,v) = \frac{e^{i(\beta - v)tn - \frac{\gamma}{2} t^2n} - 1}{i\delta t}
$$
where $\beta$ and $\gamma$ take the values $\beta_2$, $\gamma_2$ and 
$\beta_1$, $\gamma_1$, respectively.
Since the real and the imaginary part of $A_n$ are (respectively) an even and an odd function in $t$,
recalling that $\tau_n = n^{-1/2}(\log n)$ and setting $u=(\beta - v)tn$, we get
\begin{eqnarray}
	\int_{S_n} A_n(t,v) dt & = & \frac{2}{\delta}
	\int_0^{\tau_n} \frac{e^{-\frac{\gamma}{2} t^2n} \sin((\beta-v)tn)}{t}\, dt \ \ =  \nonumber \\
	\hspace{-0.3cm} = \ \frac{2}{\delta} \left\{ \int_0^{(\beta-v)\tau_n n}  \frac{\sin(u)}{u}\, du 
	\right. & - & \left.  \int_0^{(\beta-v)\tau_n n} 
	\left(1 - e^{-\frac{\gamma u^2}{2(\beta-v)^2n}}\right) \frac{\sin(u)}{u}\, du \right\}
	\label{duestar}
\end{eqnarray}
By Lemma \ref{segnialt} the first term of (\ref{duestar}) can be written as 
\begin{eqnarray}
	\nonumber
	\frac{2}{\delta} \int_0^{(\beta-v)\tau_n n} \frac{\sin(u)}{u}\, du
		& = & 
	\frac{2\mbox{ sgn}(\beta-v)}{\delta} \left( \int_0^{+\infty}\frac{\sin(u)}{u}du\; - \; 
	\int_{|\beta-v|\tau_n n}^{+\infty} \frac{\sin(u)}{u}du \right)  \\
	\label{trestar}
	& = & \ \frac{\pi}{\delta}\mbox{sgn}(\beta-v) \; - \; 
	O\left(n^{-1/2} (\log n)^{-1} \right)
\end{eqnarray}
Now we use again Lemma \ref{segnialt} to deal with the second term of (\ref{duestar}),
which has the form 
\begin{equation}
	\label{secondoadd}
	\frac{2}{\delta} \int_0^{(\beta-v)\tau_n n} B_n(u) \sin(u) du
\end{equation}
where
$B_n(u) = u^{-1} \left(1 - e^{-\frac{\gamma u^2}{2(\beta-v)^2n}}\right)$.
Note that $B_n(u)>0$ for all $u>0$, and 
$$\lim_{u\rightarrow 0} B_n(u) = 0 = \lim_{u\rightarrow +\infty} B_n(u)$$
Moreover in the set  $(0,+\infty)$ its derivative is null only at the point
$u_n = \alpha |\beta - v| \sqrt{n/\gamma}$,
for a constant $\alpha \in (1,2)$ independent of $n$ and $v$.
Thus, for $n$ large enough, $u_n$ belongs to the interval
$(0,|\beta-v| \tau_n n)$, $B_n(u)$ is increasing in the set
$(0,u_n)$ and decreasing in $(u_n,+\infty)$,
while its maximum value is
$$
B_n(u_n) = \frac{1-e^{-\frac{\alpha^2}{2}}}{\alpha|\beta-v|} \sqrt{\frac{\gamma}{n}} 
= \Theta (n^{-1/2})
$$
Defining $k_n= \lfloor\frac{u_n}{2\pi}\rfloor$ and 
$K = \lfloor \frac{|\beta - v| \tau_n n}{2\pi} \rfloor$,
we can apply Lemma \ref{segnialt} to the intervals 
$[0,2k_n\pi]$ and $[2k_n\pi+2\pi,2K\pi]$, to get
\begin{eqnarray*}
	\lefteqn{ \left| \int_0^{|\beta - v|  \tau_n n} B_n(u) \sin u\; du \right| \ \leq \ 
	2 B_n(2k_n\pi) + \left|\int_{2k_n\pi}^{2(k_n+1)\pi} B_n(u) \sin u\; du \right| \ + } \\
	& & \qquad \ + \ 2[B_n(2(k_n+1)\pi) - B_n(2K\pi)] + \int_{2K\pi}^{|\beta - v| \tau_n n} B_n(u) \sin u\; du  \ = \\
	& & \leq \ 6 B_n(u_n) \ = \ \frac{c \ \sqrt{\gamma}}{|\beta - v| \sqrt{n}} 
\end{eqnarray*}
where $c$ is a positive constant independent of $v$ and $n$.

This implies that,
for any $v$ approaching a constant different from 
$\beta_1$ and $\beta_2$,
the second term of (\ref{duestar}) is $O(n^{-1/2})$.
Therefore, applying (\ref{trestar}) and recalling that 
$v$ converges to a constant $x$ different from $\beta_1$ and $\beta_2$, we get
\begin{eqnarray*}
	\Delta_n(v) & = & \frac{2}{\delta} \ \left[ \int_{0}^{(\beta_2 -v) n \tau_n} 
	\frac{\sin u}{u} du    
	 - \int_{0}^{(\beta_1 -v) n \tau_n} 
	\frac{\sin u}{u} du \right] + \mbox{O}(n^{-1/2}) \\
& = & \frac{\pi}{\delta} \left[ \mbox{sgn}(\beta_2-v) - \mbox{sgn}(\beta_1-v) \right] + 
\mbox{O}(n^{-1/2}) \ = \ 2\pi f_U(x) + \mbox{O}(n^{-1/2})
\end{eqnarray*}
This proves equation (\ref{finale}) and hence the proof is complete.
\fdimo
The theorem clearly holds also when $k/n$ definitely lies in a finite interval
not including $\beta_1$ nor $\beta_2$ (the proof being the same).

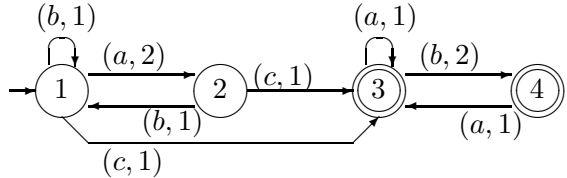
\begin{figure}[h]
	\begin{center}
		\begin{picture}(210,60)
		\put (20,30){\circle{20}}
		\put (17,27){$1$}
		\put (0,30){\vector(1,0){10}}
		\put (80,30){\circle{20}}
		\put (77,27){$2$}
		\put (140,30){\circle{20}}
		\put (137,27){$3$}
		\put (140,30){\circle{17}}
		\put (200,30){\circle{20}}
		\put (197,27){$4$}
		\put (200,30){\circle{17}}
		
		\put (20,40){\oval(10,20)[t]}
		\put (25,42){\vector(0,-1){3}}
		\put (10,55){$(b,1)$}
		\put (140,40){\oval(10,20)[t]}
		\put (145,42){\vector(0,-1){3}}
		\put (130,55){$(a,1)$}
		
		\put (30,36){\vector(1,0){40}}
		\put (35,40){$(a,2)$}
		\put (90,30){\vector(1,0){40}}
		\put (92,32){$(c,1)$}
		\put (150,36){\vector(1,0){40}}
		\put (155,40){$(b,2)$}
		
		\put (70,24){\vector(-1,0){40}}
		\put (50,15){$(b,1)$}
		\put (190,24){\vector(-1,0){40}}
		\put (170,14){$(a,1)$}
		
		\put (20,20){\line(1,-1){10}}
		\put (30,10){\line(1,0){100}}
		\put (130,10){\vector(1,1){10}}
		\put (35,0){$(c,1)$}
		\end{picture}
		\caption{Weighted finite automaton defining an equipotent bicomponent model ($\lambda_1=\lambda_2=2$) with $\beta_1 = 1/3$ and $\beta_2 = 2/3$.}
		\label{automabicomeq}
	\end{center}
\end{figure}

As an example, consider the rational stochastic model defined by the
weighted finite automaton of Figure \ref{automabicomeq},
where each transition is labelled by a pair $(\sigma,p)$, 
for a symbol $\sigma \in \{a,b,c\}$ and a weight $p>0$,
together with the arrays $\xi =(1,0,0,0)$ and $\eta=(0,0,1,1)$.
Such an automaton recognizes the set of all words $w\in\{a,b,c\}^*$
of the form $w=xcy$, such that $x,y \in \{a,b\}^*$ and
the strings $aa$ and $bb$ do not occur in $x$ and $y$, respectively.
Clearly this is a bicomponent model, with 
both pairs $(A_1,B_1)$ and $(A_2,B_2)$ aperiodic. 
Moreover $M_1=M_2$, while $A_1 \neq A_2$. 
Hence the two components are equipotent and $\beta_1 \neq \beta_2$. 
This means that $Y_n/n$ converges in distribution to a uniform r.v. of extremes
$\beta_1$, $\beta_2$, and $Y_n$ satisfies Theorem \ref{teo:bicompequipbetadiv}.
Note that simple changes may modify the limit distribution: 
for instance, setting to $3$ the weight of transition 
$2 \stackrel{b}{\rightarrow} 1$ makes dominant the first component,
implying a Gaussian local limit law (Theorem \ref{teo:locale_bicompdom}).

\subsubsection{Local limit with equal $\beta$'s and different $\gamma$'s}

In this section we present a local limit theorem for $\{Y_n\}$ defined in an equipotent communicating bicomponent model
with $\beta_1=\beta_2$ and $\gamma_1 \neq \gamma_2$.
In this case, setting $\beta=\beta_1=\beta_2$ and $\gamma= \frac{\gamma_1+\gamma_2}{2}$,
it is proved that the distribution of $\frac{Y_n - \beta n}{\sqrt{\gamma n}}$ converges 
to a mixture of Gaussian laws having  mean $0$ and variance uniformly distributed over 
an interval of extremes
$\frac{\gamma_1}{\gamma}$ and $\frac{\gamma_2}{\gamma}$ \cite{dgl04}.

Formally, we consider a r.v. $T$ having density function
\begin{equation}
\label{fconti}
f_T(x) = \frac{\gamma}{\gamma_2 - \gamma_1} 
\int_{\frac{\gamma_1}{\gamma}}^{\frac{\gamma_2}{\gamma}} \frac{e^{-\frac{x^2}{2s}}}{\sqrt{2\pi s}}\  ds
\qquad \forall\  x \in \R
\end{equation}
In passing, we observe that for each $x\in \R$, $f_T(x)$ may be regarded 
as the mean value of the ``heat kernel'' 
$\ K(x,t)=(4\pi t)^{-1/2} e^{\frac{-x^2}{4t}}\ $  at point $x$
in the time interval of extremes $\gamma_1/(2\gamma)$ and $\gamma_2/(2\gamma)$
\cite{ca84}.

Note that $E(T)=0$ and $var(T) = 1$, while its characteristic function is
\begin{equation}
\label{phicont}
\Phi_T(t) \ =\ \int_{-\infty}^{+\infty} f_T(x) e ^{itx} dx \ = \ 
2 \gamma\; \frac{e^{-\frac{\gamma_1}{2\gamma}t^2} - e^{-\frac{\gamma_2}{2\gamma}t^2}}{(\gamma_2 - \gamma_1)t^2}
\end{equation}
Clearly, $f_T(x)$ can be expressed in the form
$$
f_T(x) \ = \ \frac{1}{2\pi} \int_{-\infty}^{+\infty} \Phi_T(t) e ^{-itx} dt
\ = \ \frac{1}{2\pi} \int_{-\infty}^{+\infty}  
2 \gamma\; \frac{e^{-\frac{\gamma_1}{2\gamma}t^2} - e^{-\frac{\gamma_2}{2\gamma}t^2}}{(\gamma_2 - \gamma_1)t^2}
\ e ^{-itx} dt
$$

As in the previous section,
we assume aperiodic both pairs $(A_1,B_1)$ and $(A_2,B_2)$,
which implies Proposition \ref{terzointbic},
$c \in (0,\pi)$ is a constant 
for which relation (\ref{accanzetabicomp}) holds true and
both functions $y_1(t)$, $y_2(t)$ satisfy relations (\ref{yt}),
which can now be refined as
$$
y_j(t) \ = \ \frac{u_j(it)}{\lambda}\ = \ 
1 + i\beta t - \frac{\gamma_j + \beta^2}{2}\; t^2 + O(t^3)\ ,
\qquad \forall \ t\in \R \ : |t| \leq c, \ j=1,2
$$
Applying these values in (\ref{accanzetabicomp}), by identity (\ref{accan0bicomp})
for some $c\in(0,\pi)$ and every $t\in \R$ satisfying $0<|t| \leq c$ we obtain
\begin{equation}
\label{psingammadiv}
\Psi_n(t)\ =\ \frac{h_n(it)}{h_n(0)}\ =\  2\; \frac{1+O(t)}{1+O(1/n)}\ 
\frac{y_1(t)^n - y_2(t)^n}{(\gamma_2-\gamma_1)nt^2 + nO(t^3)} \;
+ \; \sum_{j=1,2} O\left( \frac{y_j(t)^n}{n} \right) 
\end{equation}
Now, for such a constant $c$, we split the interval $[-c,c]$ into sets $S_n$ and $V_n$ 
defined in (\ref{intervalli}).
The behaviour of $\Psi_n(t)$ in these sets is studied in the two propositions below,
where we always assume an equipotent communicating bicomponent model with 
$\beta_1 = \beta_2 = \beta$
and $\gamma_1 \neq \gamma_2$.
\begin{proposition}
	\label{secondointgammadiv}
	For some $a>0$ we have
	$|\Psi_n(t)| = o\left( n^{-3/2} \right)$ for all $t \in V_n$.
\end{proposition}
{\it Proof.}
From equation (\ref{psingammadiv}), taking $a= \min\{\gamma_1,\gamma_2\}/4$
and using (\ref{ytbound}), we can write  
\begin{equation*}
	|\Psi_n(t)| \ \leq \  O\left( \frac{|y_1(t)|^n+|y_2(t)|^n}{t^2 n} \right) \; +\;
	O\left( \sum_{j=1,2} \frac{|y_j(t)|^n}{n} \right) \ = \ 
	O\left( \frac{e^{-a(\log n)^2}}{(\log n)^2 } \right)
	\ , \  \forall\ t \in V_n
\end{equation*}
which proves the result.
\fdimo

As regards the behaviour of $\Psi_n(t)$ in $S_n$, we define
\begin{equation}
\label{secondocappaenne}
H_n(t) = 2\; \frac{e^{-\frac{\gamma_1}{2}t^2n} - e^{-\frac{\gamma_2}{2}t^2n}}{(\gamma_2-\gamma_1)t^2n} \; e^{i\beta tn}\ ,
\qquad \forall\ t \in \R
\end{equation}
It is easy to see that
${\displaystyle 
	\ |H_n(t)| \leq 2 \sum_{j=1,2} \left( 
	\frac{1 - e^{-\frac{\gamma_j}{2}t^2 n}}{|\gamma_2-\gamma_1| t^2 n} \right)\ }$
for every $t \in \R$.
Both addends take their maximum value at $t=0$, where they have a removable singularity, and such values are independent of $n$.
As a consequence we can state that 
$\ |H_n(t)| \leq \frac{\gamma_1 + \gamma_2}{|\gamma_2-\gamma_1|}$, 
for every $n\in \N_+$ and every $t\in S_n$.

\begin{proposition}
	\label{primointgammadiv}
	Let $S_n$ and $H_n(t)$ be defined by {\rm (\ref{intervalli})} and {\rm (\ref{secondocappaenne})}, respectively.
	Then, we have
	$$
	\int_{S_n} \left| \Psi_n(t) - H_n(t) \right| dt  \ = \
	O\left( n^{-1} \right)
	$$
\end{proposition}
{\it Proof.}
Starting again from equation (\ref{psingammadiv}) and applying relations 
(\ref{potenzayt}) to both $y_1(t)$ and $y_2(t)$, one can prove that
$$
\Psi_n(t)\ =\  \left[ 1 + O(t) + O(1/n) + nO(t^3) \right]\; H_n(t) \; + \; O(n^{-1} e^{-at^2 n})
$$
Applying relation (\ref{baseuno}) we get
\begin{equation}
\label{nuovo}
\int_{S_n} \left( \Psi_n(t) - H_n(t)\right)\; dt = 
 \int_{S_n}  \left(  O(t) + O(1/n) + nO(t^3) \right)H_n(t)\; dt  +  O\left( n^{-3/2} \right)
\end{equation}
Now, recalling that $H_n(t) = O(1)$ for a suitable $a>0$ we obtain the following relations
\begin{eqnarray*}
\int_{S_n} n^{-1} H_n(t)\; dt & = & O(n^{-1}) \left[ t \right]_0^{\frac{\log n}{\sqrt{n}}} = o(n^{-1}) \\
\int_{S_n} t H_n(t)\; dt & \leq & 2 \left[ t^2 \right]_0^{n^{-1/2}} + 
O \left( \int_{n^{-1/2}}^{\frac{\log n}{\sqrt{n}}} \frac{e^{-at^2n}}{\sqrt{n}}\; dt \right) 
\ = \ \Theta(n^{-1}) \\ 
\int_{S_n} nt^3 H_n(t)\; dt & = & 
\int_{0}^{\frac{\log n}{\sqrt{n}}} \Theta\left( t e^{-at^2n}\right) \; dt 
\ = \ \Theta(n^{-1})
\end{eqnarray*}
Thus, the result follows by applying the previous relations in (\ref{nuovo}).
\fdimo

We are now able to prove the main result in the present case.
\begin{theorem}
\label{teo:betaeqgammadiv}
Let $\{Y_n\}_{n\in\N}$ be defined in an equipotent communicating bicomponent model
with $\beta_1= \beta_2 = \beta$ and $\gamma_1\neq \gamma_2$.
Set $\gamma = (\gamma_1+\gamma_2)/2$ and assume aperiodic both pairs
$(A_1,B_1)$ and $(A_2,B_2)$.
Then, for $n$ tending to $+\infty$, $Y_n$ satisfies the relation
\begin{equation}
\label{convlocbugd}
\left| \sqrt{\gamma n}\; \mbox{Pr}(Y_n= k) \; - \; f_T\left(\frac{k-\beta n}{\sqrt{\gamma n}} \right) \right| \ = \ 
O\left(n^{-1/2} \right)
\end{equation}
uniformly for $k \in \{0,1,\ldots, n\}$,
where $f_T$ is defined in (\ref{fconti}).
\end{theorem}
{\it Proof.}
Again we start from equation (\ref{inversionform}) and split $[-\pi,\pi]$ into the three sets
$\{t\in\R : c < |t| \leq \pi\}$, $V_n$ and $S_n$, 
where $S_n$ and $V_n$ are defined in equalities (\ref{intervalli}), 
$c$ being a constant for which relation (\ref{psingammadiv}) holds true.
Then, by Propositions \ref{terzointbic}, \ref{secondointgammadiv} and \ref{primointgammadiv}, 
we obtain
$$
p_n(k) \ = \ \frac{1}{2\pi} \int_{S_n}
H_n(t) e^{-ikt} dt \ + \ O\left( n^{-1} \right)
$$
where $H_n(t)$ is defined in (\ref{secondocappaenne}).
Now, setting $v = \frac{k-\beta n}{\sqrt{\gamma n}}$ in the previous integral we get
\begin{equation}
\label{secondopnbic}
p_n(k) \ = \ \frac{1}{2\pi} \int_{S_n}
2 \frac{e^{-\frac{\gamma_1}{2}t^2n} - e^{-\frac{\gamma_2}{2}t^2n}}{(\gamma_2 - \gamma_1) n t^2}
e^{-iv\sqrt{\gamma n}t} dt 
\ + \ O\left( n^{-1} \right)
\end{equation}
By setting $x= t \sqrt{\gamma n}$ and recalling (\ref{phicont}), we obtain
\begin{eqnarray*}
	\lefteqn{\int_{S_n} 2
		\frac{e^{-\frac{\gamma_1}{2}t^2n} - e^{-\frac{\gamma_2}{2}t^2n}}{(\gamma_2 - \gamma_1) n t^2}
		e^{-iv\sqrt{\gamma n}t} dt  \ = \ 
	2\sqrt{\frac{\gamma}{n}}
	\int_{|x| \leq \sqrt{\gamma} \log n}
	\frac{e^{-\frac{\gamma_1}{2\gamma}x^2} - e^{-\frac{\gamma_2}{2\gamma} x^2}}{(\gamma_2 - \gamma_1) x^2}
	e^{-ixv} dx \ =} \\
	& = & 
	\frac{1}{\sqrt{\gamma n}} 
	\left\{ \int_{-\infty}^{+\infty} \Phi_T(x)e^{-ixv} dx\  - 
	\int_{|x| > \sqrt{\gamma} \log n} \Phi_T(x)e^{-ixv} dx  \right\} \  = \
	\frac{2\pi f_T(v)}{\sqrt{\gamma n}}  + o(n^{-2}) 
\end{eqnarray*}
The result follows by replacing this value in (\ref{secondopnbic}).
\fdimo

\subsubsection{Local limit with equal $\beta$'s and equal $\gamma$'s}

In this section we study the local limit properties of $\{Y_n\}$ assuming an equipotent communicating bicomponent model
with $\beta_1=\beta_2=\beta$ and $\gamma_1 = \gamma_2=\gamma$.
In this case,
it is known \cite{dgl04} that $\frac{Y_n - \beta n}{\sqrt{\gamma n}}$ converges in distribution to a
Gaussian r.v. of mean $0$ and variance $1$
and here we present a local limit law with a convergence rate of the order $O(n^{-1/2})$.

Again we assume $c\in (0,\pi)$ constant for which equality  
(\ref{svihnbicomp}) holds true, so that both functions
$y_1(t)$ and $y_2(t)$ satisfy relations (\ref{yt}) and (\ref{potenzayt}),
which can be restated as
\begin{gather}
\label{ytdij}
|y_j(t)| \leq e^{-\frac{\gamma}{4}t^2} \  \qquad  
\forall\ t\in \R \ : \ |t| \leq c, \qquad j=1,2 \\
\label{potenzaytdij}
y_j(t)^n = e^{-\frac{\gamma}{2}t^2n + i\beta t n + nO(t^3)} \  \qquad  
\forall\ t\in \R \ : \ |t| \leq n^{-q}, \qquad j=1,2
\end{gather}
where $q$ is an arbitrary value such that $1/3 < q <1/2$.

In the following two propositions the characteristic function $\Psi_n(t)$ is 
studied under conditions $|t| \leq n^{-q}$ and $ n^{-q} < |t| \leq c $, respectively,
assuming an equipotent communicating bicomponent model
with $\beta_1=\beta_2=\beta$ and $\gamma_1 = \gamma_2=\gamma$.

\begin{proposition}
	\label{secondointgammauguali} 
	For every $q\in (1/3,1/2)$, we have
	$$
	|\Psi_n(t)| =  O\left( e^{-\frac{\gamma}{4}n^{1-2q}} \right)
	\qquad \forall\ t \in \R \ : \ n^{-q} < |t| \leq c
	$$
\end{proposition}
{\it Proof.}
Applying relations (\ref{ytdij}) to equality (\ref{svihnbicomp}), we obtain
\begin{eqnarray*}
	|h_n(it)| & = & \left| s(it) \lambda^{n-1} 
	\sum_{j=0}^{n-1} y_1(t)^j y_2(t)^{n-1-j} + 
	\lambda^n \sum_{j=1,2} O(y_j(t)^n)  \right| \\
	& \leq & |s(it)|\; n\; \lambda^{n-1}\; e^{-\frac{\gamma}{4}t^2(n-1)} + 
	\lambda^n O(e^{-\frac{\gamma}{4}t^2n})
\end{eqnarray*}
and hence, by (\ref{accan0bicomp}), we have
$$
|\Psi_n(t)| \ =\ \left| \frac{h_n(it)}{h_n(0)}\right| \ \leq \  \frac{1+O(t)}{1+O(1/n)}\ 
e^{-\frac{\gamma}{4}t^2(n-1)} \; + \;  
O\left( e^{-\frac{\gamma}{4}t^2n}/n \right) 
$$
which implies the result since $n^{-q} < |t| \leq c$.
\fdimo

\begin{proposition}
	\label{primointgammauguali}
	For every $q\in (1/3,1/2)$, we have
	$$
	\int_{|t|\leq n^{-q}}
	\left| \Psi_n(t) - e^{-\frac{\gamma}{2}t^2n +i\beta tn}  \right| dt \: = \: O(n^{-1})
	$$
\end{proposition}
{\it Proof.}
From relations (\ref{svihnbicomp}) and (\ref{accan0bicomp}), applying (\ref{potenzaytdij}) and
recalling that $nO(t^3)=o(1)$ for $|t|\leq n^{-q}$, we obtain
\begin{eqnarray*}
	\Psi_n(t) \ =\  \frac{h_n(it)}{h_n(0)} & = & 
	\frac{1+O(t)}{n(1+O(1/n))}\ e^{-\frac{\gamma}{2}t^2n +i\beta tn}
	\sum_{j=0}^{n-1} e^{jO(t^3)+(n-1-j)O(t^3)} \; + \; O\left( \frac{e^{-\frac{\gamma}{2}t^2n}}{n} \right) \\
	& = & \left(1 + O(t) + O(n^{-1}) +nO(t^3) \right)e^{-\frac{\gamma}{2}t^2n +i\beta tn} \; + \;
	O\left( e^{-\frac{\gamma}{2}t^2n} / n \right)
\end{eqnarray*}
Therefore, a straightforward computation shows that
\begin{eqnarray*}
\lefteqn{
\int_{|t|\leq n^{-q}}
\left| \Psi_n(t) - e^{-\frac{\gamma}{2}t^2n +i\beta tn}  \right| dt \; = \;
} \\ \qquad & = & 
\int_0^{n^{-q}} \left(O(t)+O(n^{-1})+nO(t^3)\right) e^{-\frac{\gamma}{2}t^2n} dt + O(n^{-1-q})
 \; = \;  O(n^{-1})
\end{eqnarray*}
\fdimo

Now we are able to state the local limit theorem in the present case.
For the proof one can argue as in Theorem \ref{teo:locale_primitivo},
replacing Propositions \ref{terzoint}, \ref{secondoint} and \ref{intpsint} 
by Propositions
\ref{terzointbic}, \ref{secondointgammauguali} and \ref{primointgammauguali}, 
respectively.

\begin{theorem}
	\label{teo:loclimgammauguali}
	Let $\{Y_n\}_{n\in\N}$ be defined in an equipotent communicating bicomponent model
	with $\beta_1= \beta_2 = \beta$ and $\gamma_1 = \gamma_2 = \gamma$, 
	and assume aperiodic both pairs $(A_1,B_1)$ and $(A_2,B_2)$.
	Then, for $n$ tending to $+\infty$ the relation
	$$
	\left| \sqrt{n} \mbox{Pr}\left(Y_n = k \right) \: - \: 
	\frac{e^{-\frac{(k-\beta n)^2}{2\gamma n}}}{\sqrt{2\pi\gamma}} \right|
	\: = \: \mbox{O}\left(n^{-1/2}\right)
	$$
	holds true uniformly for every $k \in \{0,1,\ldots,n\}$.
\end{theorem}

\section{Sum models}
\label{sec:sum}

In this section we study the problem assuming a bicomponent rational model without communication.
Formally, the linear representation $(\xi,A,B,\eta)$ defining $\{Y_n\}$ satisfies conditions 1. and 2. of
Section \ref{sec:bicomponente}, for two suitable 4-tuples
$(\xi_1,A_1,B_1,\eta_1)$, $(\xi_2,A_2,B_2,\eta_2)$, together with the further condition
\begin{description}
	\item[3'.] $\xi_1 \neq 0 \neq \eta_1$, $\xi_2 \neq 0 \neq \eta_2$ and $A_0= [0] = B_0$.
\end{description}
In this case, for every $w\in\{a,b\}^*$ we have
$$
\xi' \mu(w) \eta = \xi'_1 \mu_1(w) \eta_1 + \xi'_1 \mu_1(w) \eta_1
$$
where $\mu$, $\mu_1$ and $\mu_2$ are the morphisms defined by pairs
$(A,B)$, $(A_1,B_1)$ and $(A_2,B_2)$, respectively.
This means that the formal series $r$ with linear representation $(\xi,A,B,\eta)$ is 
the sum of two rational formal series $r_1$, $r_2$ with irreducible linear representation,
i.e. $(r,w) = (r_1,w)+(r_2,w)$ for every $w\in\{a,b\}^*$.

Under these hypotheses, for sake of brevity, we say that $\{Y_n\}_n$ is defined in a \emph{sum model}.
Adopting the same notation of Section \ref{sec:bicomponente}, 
here we have $G(x,y)=0$ in relations (\ref{H1eH2}) implying, 
for every $z \in \C$ and $t \in \R$, the identities
\begin{equation}
\label{accaennesum} 
h_n(z) = h_n^{(1)}(z) + h_n^{(2)}(z)
\qquad  \Psi_n(it) = \frac{h_n^{(1)}(it) +h_n^{(2)}(it)}{h_n(0)}
\end{equation}

Again the simplest case occurs when there exists a dominant component.
Recall that in this case $\{Y_n\}$ has a Gaussian limit distribution \cite{dgl04}
and this result can be extended to a local limit law as stated in the following statement,
whose proof is similar to that one of Theorem \ref{teo:locale_bicompdom}.

\begin{theorem}
	\label{teo:locale_binocomdom}
	Let $\{Y_n\}$ be defined in a sum model with
	$\lambda_1 > \lambda_2$ and $M_1$ aperiodic (and hence primitive).
	Also assume aperiodic the pair $(A_1,B_1)$.
	Then $0<\beta_1 < 1$, $0<\gamma_1$ and, as $n$ tends to $+\infty$, the relation 
	$$
	\left| \sqrt{n} \mbox{Pr}\left(Y_n = k \right) \: - \: 
	\frac{e^{-\frac{(k-\beta_1 n)^2}{2\gamma_1 n}}}{\sqrt{2\pi\gamma_1}} \right|
	\: = \: \mbox{O}\left(n^{-1/2}\right)
	$$
	holds true uniformly for every $k\in \{0,1,\ldots,n\}$.
\end{theorem}

 \subsection{Equipotent sum models}
 
 \label{sec:equipotentsm}
 We now study the local limit properties of our statistics for non-communi\-ca\-ting 
 bicomponent models in the equipotent case.
 More precisely, let $\{Y_n\}$ be defined in a sum model with $\lambda_1=\lambda_2 = \lambda$ and 
 both matrices $M_1$, $M_2$ aperiodic (and hence primitive).
 Under these hypotheses we say that $\{Y_n\}$ is defined in an \emph{equipotent sum model}.
 The limit distribution of $\{Y_n\}$ in this case is studied in \cite{dgl04} and
 depends on the parameters $\alpha_1$, $\alpha_2$, $\beta_1$, $\beta_2$, $\gamma_1$, $\gamma_2$
 defined in (\ref{betaegamma}).
 Here we prove local limit properties, with a convergence rate $O(n^{-1/2})$,
 under the further assumption that both pairs $(A_1,B_1)$ and $(A_2,B_2)$ are aperiodic.
 To this end we first determine some identities for function $h_n(z)$ in the present case.

 By properties of the primitive matrices \cite{se81} it is easy to see that
 \begin{eqnarray*}
 	h_n(0) = \xi' M^n \eta & = & \xi'_1 \nu_1 \zeta_1' \eta_1 \cdot \lambda^n + 
 	\xi'_2 \nu_2 \zeta_2'\eta_2 \cdot \lambda^n + O(\rho^n) \\
 	& = & (\alpha_1 + \alpha_2) \lambda^n + O(\rho^n) \ ,
 	\qquad 0\leq \rho < \lambda 
 \end{eqnarray*}
 where $\zeta_j$ and $\nu_j$ are the eigenvectors defined in Section \ref{sec:primitive},
 for $j=1,2$.
 Also note that $\alpha_j = r_j(0)$ for each $j$, $r_j(z)$ being the same as in (\ref{sviluppo-h-eq}).
 Using these identities function $\Psi_n(t)$ can be evaluated from
 (\ref{accaennesum}).
 
Clearly, also the type of local limit law we present depends on parameters
 $\alpha_1$, $\alpha_2$,
 $\beta_1$, $\beta_2$, $\gamma_1$, $\gamma_2$.
 In general we obtain a local limit law towards a convex combination of two Gaussian distributions,
 which coincide when $\beta_1 = \beta_2$ and $\gamma_1 = \gamma_2$.
 
 \begin{theorem}
 	\label{teo:locale_bicompequi}
 	Let $\{Y_n\}$ be defined in an equipotent sum model
 	and assume that both pairs $(A_1,B_1)$, $(A_2,B_2)$ are aperiodic.
 	Then, as $n$ tends to $+\infty$, the relation 
 	$$
 	\left| \sqrt{n} \mbox{Pr}\left(Y_n = k \right) \: - \: 
 	\left( \frac{\alpha_1}{\alpha_1+\alpha_2} \frac{e^{-\frac{(k-\beta_1 n)^2}{2\gamma_1 n}}}{\sqrt{2\pi\gamma_1}} +
 	\frac{\alpha_2}{\alpha_1+\alpha_2} \frac{e^{-\frac{(k-\beta_2 n)^2}{2\gamma_2 n}}}{\sqrt{2\pi\gamma_2}}\right) 
 	\right| \: = \: \mbox{O}\left(n^{-1/2}\right)
 	$$
 	holds true uniformly for every $k\in \{0,1,\ldots,n\}$.
 \end{theorem}
 {\bf Proof.}
 Again the main idea is to study the characteristic function $\Psi_n(t)$ for
 $t \in [-\pi, \pi]$ by splitting this interval into the three sets given in 
 (\ref{intervalpr}),
 where $c\in (0,\pi)$ is a constant satisfying relations (\ref{yt})
 for both $y_1(t)$ and $y_2(t)$,
 and  $q$ is an arbitrary value such that $\frac{1}{3} < q < \frac{1}{2}$.
 The behaviour of $\Psi_n(t)$ in these sets is characterized by the following properties:
 
 {\bf a.}
 For some $\varepsilon \in(0, 1)$ we have
 \begin{equation}
 \label{terzointbis}
 |\Psi_n(t)| = O(\varepsilon^n)\  \qquad  
 \forall\ t\in \R \ : \  c < |t| \leq \pi
 \end{equation}
 
 {\bf b.}
 There exists $a>0$ such that
 \begin{equation}
 \label{secondointbis}
 |\Psi_n(t)| = O\left( e^{- a n^{1-2q}} \right)
 \qquad \forall\ t\in \R \ : \  n^{-q} < |t| \leq c
 \end{equation}
 
 {\bf c.}
 \begin{equation}
 \label{intpsintbis}
 \int_{|t|\leq n^{-q}}
 	\left| \Psi_n(t) - \frac{\alpha_1}{\alpha_1+\alpha_2} e^{-\frac{\gamma_1}{2}t^2n +i\beta_1 tn} - \frac{\alpha_2}{\alpha_1+\alpha_2} e^{-\frac{\gamma_2}{2}t^2n +i\beta_2 tn} \right| dt \: = \: O(n^{-1})
 \end{equation}
 
 \noindent
 {\it Proof of }(\ref{terzointbis}).
 We can argue as in Proposition \ref{terzointbicomp}.
 The only difference is that now the eigenvalues of $A_2e^{it} + B_2$ are smaller than
 $\lambda = \lambda_2$, and this simplifies the proof.
 
 \smallskip
 \noindent
 {\it Proof of }(\ref{secondointbis}).
 By relation (\ref{sviluppo-h-eq}),  
 for some $\varepsilon \in (0,1)$ and all $t\in \R$ satisfying $|t| \leq c$, we have
 \begin{equation}
 \label{psiditinzerobis}
 \Psi_n(t) = \frac{h_n(it)}{h_n(0)} = 
 \frac{r_1(it) u_1(it)^n + r_2(t)u_2(it)^n}{(r_1(0) + r_2(0)) \lambda^n} +
 \mbox{\small $O(\varepsilon^n)$}
 = \sum_{j=1,2} c_j y_j(t)^n  + \mbox{\small $O(\varepsilon^n)$}
 \end{equation}
 where $c_1$ and $c_2$ are positive constants. 
 Also, setting $a= \min\{\gamma_1/4, \gamma_2/4 \}$, by inequality (\ref{yt}) recalling
 $n^{-q} \leq |t| \leq c$ we obtain  
 $|y_j(t)|^n \leq e^{-a n^{1-2q}}$, for each $j=1,2$,
 which implies the result.
 
 \smallskip
 \noindent
 {\it Proof of }(\ref{intpsintbis}).
 From equality (\ref{psiditinzerobis}),
 applying relation (\ref{potenzayt}) and recalling that $nO(t^3) = o(1)$ for $|t| \leq n^{-q}$, in the same interval for $t$ we get
 $$
 \Psi_n(t) = \sum_{j=1,2}\frac{r_j(0)+ O(t)}{r_1(0)+r_2(0)}(1 + nO(t^3)) e^{-\frac{\gamma_j}{2}t^2n +i\beta_j tn} +
 O(\varepsilon^n)
 $$
 Thus, since $r_j(0) = \alpha_j$ for each $j$, reasoning as in the proof of Proposition \ref{intpsint}
 we obtain
 \begin{eqnarray}
 \nonumber
 {\lefteqn{\int_{|t|\leq n^{-q}}
 		\left| \Psi_n(t) -  \sum_{j=1,2} 
 		\frac{\alpha_j}{\alpha_1+\alpha_2} e^{-\frac{\gamma_j}{2}t^2n +i\beta_j tn} \right| dt =} }
 \\ 
 & = &  \sum_{j=1,2}
 \int_{|t|\leq n^{-q}} |O(t) + nO(t^3)|\; e^{-\frac{\gamma_j}{2}t^2n} dt +
 O(\varepsilon^n)  = O\left( n^{-1} \right) 
 \end{eqnarray}
 Now consider our main goal. 
 Defining $p_n(k)= \mbox{Pr}\left\{Y_n = k \right\}$,
 from the inversion formula (\ref{inversionform}), 
 by relations (\ref{terzointbis}), (\ref{secondointbis}) and (\ref{intpsintbis}),
 we obtain
 \begin{eqnarray}
 \nonumber
 p_n(k) & = & \frac{1}{2\pi} \int_{|t|\leq n^{-q}} \Psi_n(t) e^{-itk} dt +
 O\left( e^{-a n^{1-2q}} \right) + 
 O(\varepsilon^n) \\ & = &
 \label{pnkintbis}
 \frac{1}{2\pi} \sum_{j=1,2} \frac{\alpha_j}{\alpha_1+\alpha_2}
 \int_{|t|\leq n^{-q}} e^{-\frac{\gamma_j}{2}t^2n +i\beta_j tn - itk} dt + O(n^{-1})
 \end{eqnarray}
 Moreover, defining the variables $v_j = \frac{k-\beta_j n}{\sqrt{\gamma_j n}}$,
 for $j=1,2$,
 the last integrals can be evaluated as in (\ref{due}) and (\ref{uno}), obtaining
 $$
 \int_{|t|\leq n^{-q}} e^{-\frac{\gamma_j}{2}t^2n +i\beta_j tn - itk} dt
 = 
 \frac{1}{\sqrt{\gamma_j n}} \left( \sqrt{2\pi}\  e^{-\frac{v_j^2}{2}}
 + O(e^{-\frac{\gamma_j}{2}n^{1-2q} } ) \right)
 $$
 which replaced in (\ref{pnkintbis}) yields the result.
 \fdimo
 
 We observe that if $\beta_1 = \beta_2$ and $\gamma_1 = \gamma_2$
 then the limit density given by Theorem \ref{teo:locale_bicompequi} 
 reduces to a Gaussian law.
 This yields the following
 
 \begin{corollary}
 	\label{cor:equipgaus}
 	Let $\{Y_n\}$ be defined in an equipotent sum model
 	with $\beta_1 = \beta_2 = \beta$, $\gamma_1 = \gamma_2 = \gamma$
 	and assume aperiodic both pairs $(A_1,B_1)$, $(A_2,B_2)$.
 	Then, as $n$ tends to $+\infty$, the relation 
 	$$
 	\left| \sqrt{n} \mbox{Pr}\left(Y_n = k \right) \: - \: 
 	\frac{e^{-\frac{(k-\beta n)^2}{2\gamma n}}}{\sqrt{2\pi\gamma}} \right|
 	\: = \: \mbox{O}\left(n^{-1/2}\right)
 	$$
 	holds true uniformly for every $k\in \{0,1,\ldots,n\}$.
 \end{corollary}
 
 On the contrary, when $\beta_1 \neq \beta_2$ or $\gamma_1 \neq \gamma_2$ (or both)
 the previous result yields a local limit law toward a convex 
 combination of two Gaussian distributions that differ by 
 their mean value or by their variance. 
 More precisely, in this case we obtain the distribution of a r.v. $\cal L$ defined by
 \begin{equation}
 \label{legge_elle}
 {\cal L} = \left[\beta_1 {\cal B}_p + \beta_2 ( 1 - {\cal B}_p )\right] n +
 \left[{\cal B}_p {\cal N}_{0,\gamma_1 } + (1 - {\cal B}_p ){\cal N}_{0,\gamma_2 }\right] \sqrt{n}
 \end{equation}
 where ${\cal B}_p$ is a Bernoullian r.v. of parameter 
 $p =  \frac{\alpha_1}{\alpha_1+\alpha_2}$, and
 ${\cal N}_{0,\gamma_j}$ is a Gaussian r.v. of mean $0$ and variance $\gamma_j$,
 assuming mutually independent all variables ${\cal B}_p$, ${\cal N}_{0,\gamma_1}$, ${\cal N}_{0,\gamma_2}$.
 In particular from the analysis above it is clear that
 $$
 \frac{Y_n - \left[\beta_1 {\cal B}_p + \beta_2 ( 1 - {\cal B}_p )\right] n}{[{\cal B}_p\sqrt{\gamma_1} + (1-{\cal B}_p)\sqrt{\gamma_2}]\sqrt{n}}\ \longrightarrow {\cal N}_{0,1} \qquad
 \mbox{ in distribution} 
 $$
 which specializes into the following cases:
 \begin{description}
 \item[-] if $\beta_1 = \beta_2 = \beta$ and $\gamma_1 \neq \gamma_2$ then 
 $$
 \frac{Y_n - \beta n}{\sqrt{n}}\ \longrightarrow 
 [{\cal B}_p{\cal N}_{0,\gamma_1} + (1-{\cal B}_p){\cal N}_{0,\gamma_2}] \qquad
 \mbox{ in distribution} 
 $$
 \item[-] if $\beta_1 \neq \beta_2$ and $\gamma_1 = \gamma_2 = \gamma$ then
 $$
 \frac{Y_n - \left[\beta_1 {\cal B}_p + \beta_2 ( 1 - {\cal B}_p )\right] n}{\sqrt{n}}\ \longrightarrow 
 {\cal N}_{0,\gamma} \qquad
 \mbox{ in distribution} 
 $$
\end{description}
 A curious fact is that $\cal L$ as defined in (\ref{legge_elle}) also depends on the weights of 
 initial and final states ($\xi$, $\eta$).
 This does not occur in the equipotent bicomponent models with communication, studied in Section \ref{subsec:equipotent},
 and seems to suggest that the present model is not ergodic
 (in the sense that the limit distribution depends on the starting states).
 
\begin{figure}[h]
 	\begin{center}
 		\begin{picture}(210,60)
 		\put (20,30){\circle{20}}
 		\put (17,27){$1$}
 		\put (0,30){\vector(1,0){10}}
 		\put (80,30){\circle{20}}
 		\put (77,27){$2$}
 		\put (80,30){\circle{17}}
 		\put (140,30){\circle{20}}
 		\put (137,27){$3$}
 		\put (140,30){\circle{17}}
 		\put (120,30){\vector(1,0){10}}
 		\put (200,30){\circle{20}}
 		\put (197,27){$4$}
 		\put (200,30){\circle{17}}
 		\put (20,40){\oval(10,20)[t]}
 		\put (25,42){\vector(0,-1){3}}
 		\put (10,55){$(b,1)$}
 		\put (140,40){\oval(10,20)[t]}
 		\put (145,42){\vector(0,-1){3}}
 		\put (130,55){$(a,1)$}
 		\put (30,36){\vector(1,0){40}}
 		\put (35,40){$(a,2)$}
 		\put (150,36){\vector(1,0){40}}
 		\put (155,40){$(b,2)$}
 		\put (70,24){\vector(-1,0){40}}
 		\put (50,15){$(b,1)$}
 		\put (190,24){\vector(-1,0){40}}
 		\put (170,14){$(a,1)$}
 		\end{picture}
 		\vspace{-0.5cm}
 		\caption{Weighted finite automaton defining a non-communicating bicomponent model with $\lambda_1=\lambda_2=2$, $\alpha_1=2/3$, $\alpha_2 =4/3$, 
 			$\beta_1=1/3$, $\beta_2 = 2/3$, $\gamma_1 = \gamma_2 = 2/27$.}
 		\label{automabinocomeq}
 	\end{center}
 \end{figure}
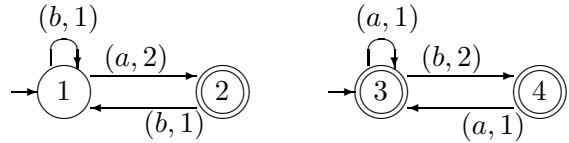

 As an example, consider the rational model defined by the
 weighted finite automaton of Figure \ref{automabinocomeq},
 together with $\xi =(1,0,1,0)$ and $\eta=(0,1,1,1)$.
 Such an automaton recognizes the set of all words $\{w\in\{a,b\}^*$ 
 such that pattern $aa$ or pattern $bb$ (or both) do not occur in $w$.
 Clearly this is a bicomponent model, with 
 both pairs $(A_1,B_1)$ and $(A_2,B_2)$ aperiodic. 
 Moreover $M_1=M_2$, while $A_1 \neq A_2$. 
 Hence the two components are equipotent and $\beta_1 \neq \beta_2$,
 however one can show that $\gamma_1 = \gamma_2$ \cite{bcgl03}.
 This implies a local limit law towards a convex combination of 
 two Gaussian laws having different main constant of the mean value 
 (but equal main constant of the variance).
 Note that simple changes may modify the limit distribution: 
 for instance, setting to $2$ the weight of transition 
 $2 \stackrel{b}{\rightarrow} 1$ makes dominant the first component,
 implying a Gaussian local limit law (Theorem \ref{teo:locale_binocomdom}).
 
\section{Conclusions}
\label{sec:conclusioni}

In this work we have studied the local limit laws of symbol statistics defined in primitive rational models and in bicomponent rational models with or without communication.
These laws, summarized in Table \ref{tavola}, yield a convergence rate of the order $O(n^{-1/2})$ and are obtained assuming suitable aperiodicity conditions concerning the number of symbol occurrences in cycles of equal length.

\begin{table}[h]
	\label{tavola}
	\begin{center}
		{
			\begin{tabular}{|l|l|l|l|l|l|}
				\hline
				& {$\begin{array}{c} \mbox{ } \\ \mbox{Primitive} \\ \mbox{models} \end{array} $} & \multicolumn{4}{c|}{$\begin{array}{c} \mbox{Bicomponent models} \\ \mbox{with communication} \end{array} $  } \\ \cline{3-6} 
				{\bf 1.} &  & dominant & \multicolumn{3}{c|}{equipotent} \\ \cline{4-6}
				&  &  & \multicolumn{1}{c|}{$\beta_1\neq\beta_2$} & 
				\multicolumn{1}{c|}{$\begin{array}{c} \beta_1=\beta_2 \\ \gamma_1\neq\gamma_2 \end{array}$} &
				$\begin{array}{c} \beta_1=\beta_2 \\ \gamma_1=\gamma_2 \end{array}$ \\
				\hline \hline 
				$\begin{array}{c} \mbox{Local limit} \\ \mbox{distribution} \end{array}$
				& ${\cal N}_{0,1}$  & ${\cal N}_{0,1}$ & $U_{\beta_1,\beta_2}$ & 
				$\ \ T$ & ${\cal N}_{0,1}$ \\
				\hline \hline
				&  \multicolumn{5}{c|}{$\begin{array}{c} \mbox{Bicomponent models} \\ \mbox{without communication} \end{array} $  } \\ \cline{2-6}
				{\bf 2.} & dominant  & \multicolumn{4}{c|}{equipotent} \\ \cline{3-6}
				& & \multicolumn{1}{c|}{$\begin{array}{c} \beta_1\neq\beta_2 \\ \gamma_1\neq\gamma_2 \end{array}$} 
				& \multicolumn{1}{c|}{$\begin{array}{c} \beta_1=\beta_2 = \beta \\ \gamma_1\neq\gamma_2 \end{array}$} & 
				\multicolumn{1}{c|}{$\begin{array}{c} \beta_1\neq\beta_2 \\ \gamma_1=\gamma_2 = \gamma \end{array}$} &
				$\begin{array}{c} \beta_1=\beta_2 \\ \gamma_1=\gamma_2 \end{array}$ \\
				\hline \hline 
				$\begin{array}{c} \mbox{Local limit} \\ \mbox{distribution} \end{array}$
				& ${\cal N}_{0,1}$ & ${\cal L}$ & {\small ${\cal B}_p {\cal N}_{0,\gamma_1 } + (1 - {\cal B}_p ){\cal N}_{0,\gamma_2 }$} & 
				${\cal N}_{0,\gamma}$ & ${\cal N}_{0,1}$ \\
				\hline
		\end{tabular}} \medskip
		\caption{Symbols ${\cal N}_{0,1}$, $U_{\beta_1,\beta_2}$ and $T$ denote respectively
			a Gaussian, uniform and $T$-type local limit, 
			$T$ being defined in Section 4.2.2.
			Also, the r.v.'s ${\cal L}$ and ${\cal B}_p$ are defined in
		(\ref{legge_elle}) .} 
	\end{center}
\end{table}

Our analysis of the bicomponent models includes the main cases but is not exhaustive,
since it does not contain the degenerate cases, i.e.
when either $A_i=0$ or $B_i=0$ for a dominant component $i\in\{1,2\}$.
In these cases a large variety of different limit distributions is obtained
\cite[Section 8]{dgl04}, most of them related to natural matrix extensions of the geometric distribution.
This large range of possible limit behaviours is also mentioned in \cite{bbpt12} and it seems natural to
study common properties of these distributions, determining a suitable classification.

\bibliographystyle{plain}


\end{document}